\documentclass[11pt]{article} \makeatother
\usepackage{amsfonts,amsmath,amssymb}
\usepackage{graphicx}
\usepackage{mathrsfs}
\allowdisplaybreaks[4]
\newcommand{\lbl}[1]{\label{#1}}
\newtheorem{theo}{Theorem}[section]
\newtheorem{lemma}[theo]{Lemma}

\newcommand{\be}{\begin{equation}}
\newcommand{\ee}{\end{equation}}
\newcommand\bes{\begin{eqnarray}} \newcommand\ees{\end{eqnarray}}
\newcommand{\bess}{\begin{eqnarray*}}
\newcommand{\eess}{\end{eqnarray*}}

\numberwithin{equation}{section}

\begin{document}

 \begin{center} {\bf\Large Existence and general stabilization of the Timoshenko system }\\[2mm]
{\bf\Large with a thermo-viscoelastic damping and a delay term in the}\\[2mm]
{\bf\Large internal feedback}\\[2mm]

{\large Weican Zhou and Miaomiao Chen \footnote{Corresponding
    author. E-mail: {\sf mmchennuist@163.com}.} \\[1mm]}
{College of Mathematics and Statistics, Nanjing University of
Information Science and Technology, Nanjing 210044, China}\\[2mm]
\end{center}

\setlength{\baselineskip}{17pt}{\setlength\arraycolsep{2pt}

\begin{quote}
\noindent {\bf Abstract:} In this paper, we consider a Timoshenko system with a thermo-viscoelastic damping and a delay term in the internal feedback
$$ \left\{
\begin{array}{lll}\displaystyle
\rho_1\varphi_{tt}-K(\varphi_{x}+\psi)_{x}=0, \ \ &(x,t)\in (0,1)\times(0,\infty), \medskip\\
\rho_2\psi_{tt}-b\psi_{xx}+K(\varphi_{x}+\psi)+\beta\theta_{x}=0, \ \ &(x,t)\in (0,1)\times(0,\infty), \medskip\\
\displaystyle \rho_3\theta_{tt}-\delta\theta_{xx}+\gamma\psi_{ttx}+\int_{0}^{t}g(t-s)\theta_{xx}(s){\rm d}s\medskip\\
\quad\quad+\mu_{1}\theta_{t}(x,t)+\mu_{2}\theta_{t}(x,t-\tau)=0, \ \ &(x,t)\in (0,1)\times(0,\infty)
\end{array}
\right.\lbl{0.1}
$$
together with initial datum and boundary conditions of Dirichlet type, where $g $ is a positive non-increasing relaxation function and $\mu_1, \mu_2 $ are positive constants. Under an hypothesis between the weight of the delay term in the feedback and the the weight of the friction damping term, using the Faedo-Galerkin approximations together with some energy estimates, we prove the global existence of the solutions. Then, by introducing appropriate Lyapunov functionals, under the imposed constrain on the weights of the two feedbacks and the coefficients, we establish the general energy decay result from which the exponential and polynomial types of decay are only special cases.

\noindent {\bf Keywords}: {Global existence; General energy decay; Timoshenko system; Relaxation function.}

\noindent {\bf AMS Subject Classification (2010):} {\small 35B35,
35L55,  74D05, 93D15}
\end{quote}

\setlength{\baselineskip}{17pt}{\setlength\arraycolsep{2pt}

\section{Introduction }\label{s1}
In this paper we investigate the existence and decay properties of solutions for the Timoshenko system with a thermo-viscoelastic damping and a delay term of the form
\be \left\{
\begin{array}{lll}\displaystyle
\rho_1\varphi_{tt}-K(\varphi_{x}+\psi)_{x}=0, \ \ &(x,t)\in (0,1)\times(0,\infty), \medskip\\
\rho_2\psi_{tt}-b\psi_{xx}+K(\varphi_{x}+\psi)+\beta\theta_{x}=0, \ \ &(x,t)\in (0,1)\times(0,\infty), \medskip\\
\displaystyle\rho_3\theta_{tt}-\delta\theta_{xx}+\gamma\psi_{ttx}+\int_{0}^{t}g(t-s)\theta_{xx}(s){\rm d}s\medskip\\
 \quad\quad+\mu_{1}\theta_{t}(x,t)+\mu_{2}\theta_{t}(x,t-\tau)=0, \ \ &(x,t)\in (0,1)\times(0,\infty)
\end{array}
\right.
\lbl{1.1}
\ee
with the following initial datum and boundary conditions:
\be \left\{
\begin{array}{lll}\displaystyle
\varphi(x,0)=\varphi_0, \ \ \varphi_{t}(x,0)=\varphi_1, \ \ \psi(x,0)=\psi_0, \ \ \psi_{t}(x,0)=\psi_1,\ \ &x\in [0,1], \medskip\\
\theta(x,0)=\theta_0, \ \ \theta_{t}(x,0)=\theta_1,\ \ &x\in [0,1], \medskip\\
\varphi(0,t)=\varphi(1,t)=\psi(0,t)=\psi(1,t)=\theta_{x}(0,t)=\theta_{x}(1,t)=0,\ \ &t\in [0,\infty), \medskip\\
\theta_{t}(x,t-\tau)=f_{0}(x,t-\tau), \ \ & (x,t)\in (0,1)\times(0,\tau),
\end{array}
\right.
\lbl{1.11}
\ee
where the coefficients  $\rho_1, \rho_2, \rho_3, K, b, \beta, \gamma, \delta, \mu_{1} $ and $\mu_2 $ are positive constants, $\tau>0 $ represents the time delay.

System \eqref{1.1} arises in the theory of the transverse vibration of a beam which was first introduced by Timoshenko. In 1921, Timoshenko \cite{timo1921} considered the following system of coupled hyperbolic equations:
\be \left\{
\begin{array}{lll}\displaystyle
\rho \varphi_{tt}-K\left(\varphi_{x}-\psi\right)_{x}=0, \ \ & (x,t)\in (0,L)\times(0,\infty), \medskip\\
I_{\rho}\psi_{tt}-\left(EI\psi_{x}\right)_{x}-K(\varphi_{x}-\psi)=0, \ \ & (x,t)\in (0,L)\times(0,\infty),
\end{array}
\right.\lbl{10.5}
\ee
where $\varphi $ is the transverse displacement of the beam and $\psi $ is the rotation angle of the filament of the beam. The coefficients $\rho, I_\rho, E, I $ and $K $ are the density, the polar moment of inertia of a cross section, Young's modulus of elasticity, the moment of inertia of a cross section, and the shear modulus, respectively.

Many mathematicians have studied system \eqref{10.5} and some results concerning the existence and asymptotic behavior of solutions have been established, see for instance \cite{djebablat2010,liu1998,shi1998,shubov2002} and the references therein. Kim and Renardy \cite{kim1987} considered \eqref{10.5} together with two linear boundary conditions of the form
\bess
K\psi(L,t)-K\frac{\partial\varphi}{\partial x}(L,t)=\alpha\frac{\partial\varphi}{\partial t}(L,t),\ \ &t\in [0,\infty), \medskip\\
EI\frac{\partial\psi}{\partial x}(L,t)=-\beta\frac{\partial\psi}{\partial t}(L,t),\ \ &t\in [0,\infty),
\eess
and used the multiplier techniques to establish an exponential decay result for the energy of \eqref{10.5}. They also provided numerical estimates to the eigenvalues of the operator associated with the system \eqref{10.5}. Soufyane and Wehbe \cite{soufyane2003} considered
\be \left\{
\begin{array}{lll}\displaystyle
\rho \varphi_{tt}-K\left(\varphi_{x}-\psi\right)_{x}=0, \ \ & (x,t)\in (0,L)\times(0,\infty), \medskip\\
I_{\rho}\psi_{tt}-\left(EI\psi_{x}\right)_{x}-K(\varphi_{x}-\psi)+b(x)\psi_{t}=0, \ \ & (x,t)\in (0,L)\times(0,\infty),
\end{array}
\right.\lbl{0.00}
\ee
where $b $ is a positive and continuous function, satisfying
$$b(x)\geq b_0>0, \ \  \forall x\in [a_0,a_1]\subset [0,L].$$
They proved that the uniform stability of \eqref{0.00} holds if and only if the wave speeds are equal $\left(\frac{K}{\rho}=\frac{EI}{I_{\rho}}\right); $ Otherwise only the asymptotic stability has been proved. For more results related to system \eqref{10.5}, we refer the readers to \cite{feng1998,raposo2005,santos2002,riverra2003} and references therein.

In \cite{messaoudi2008}, Messaoudi and Said-Houari considered the  one-dimensional linear Timoshenko system of thermoelastic type
\bes \left\{
\begin{array}{lll}\displaystyle
\rho_1\varphi_{tt}-K(\varphi_{x}+\psi)_{x}=0, \ \ & (x,t)\in (0,1)\times(0,\infty), \medskip\\
\rho_2\psi_{tt}-b\psi_{xx}+K(\varphi_{x}+\psi)+\beta\theta_{x}=0, \ \ & (x,t)\in (0,1)\times(0,\infty), \medskip\\
\rho_3\theta_{tt}-\delta\theta_{xx}+\gamma\psi_{ttx}-\kappa\theta_{txx}=0, \ \ & (x,t)\in (0,1)\times(0,\infty).
\end{array}
\right.\lbl{0.0}
\ees
They used the energy method to prove an exponential decay  under the condition $\frac{\rho_1}{K}=\frac{\rho_2}{b}. $  A similar result was also  obtained  by Rivera and Racke \cite{riverra2002} and Messaoudi et al. \cite{messaoudipo2009}. Then, in \cite{messaoudif2013}, Messaoudi and Fareh also considered problem \eqref{0.0}. By introducing the first and second-order energy functions, they proved a polynomial stability result under the condition $\frac{\rho_1}{K}\neq\frac{\rho_2}{b}. $

The case of time delay in the Timoshenko system has been studied by some authors. Said-Houari and Laskri \cite{houari2010} considered the following Timoshenko system with a constant time delay in the feedback:
\be
\left\{
\begin{array}{lll}\displaystyle
\rho_{1}\varphi_{tt}-K\left(\varphi_{x}+\psi\right)_{x}=0, \ \ & (x,t)\in (0,1)\times(0,\infty), \medskip\\
\rho_{2}\psi_{tt}-b\psi_{xx}+K\left(\varphi_{x}+\psi\right)+\mu_1\psi_{t}+\mu_{2}\psi_{t}(x,t-\tau)=0, \ \ & (x,t)\in  (0,1)\times(0,\infty).
\end{array}
\right.
\lbl{1.2}
\ee
They established an exponential decay result for the case of equal-speed wave propagation $\left(\frac{\rho_1}{K}=\frac{\rho_2}{b}\right) $ under the assumption $\mu_2<\mu_1.$ Then, Kirane el al. \cite{houaria2011} considered the Timoshenko system with a time-varying delay
$$ \left\{
\begin{array}{lll}\displaystyle
\rho_1\varphi_{tt}-K\left(\varphi_{x}+\psi\right)_{x}=0, \ \ & (x,t)\in (0,1)\times(0,\infty), \medskip\\
\rho_2\psi_{tt}-b\psi_{xx}+K\left(\varphi_{x}+\psi\right)+\mu_1\psi_t+\mu_2\psi_{t}\left(x,t-\tau(t)\right)=0, \ \ & (x,t)\in  (0,1)\times(0,\infty),
\end{array}
\right.
$$
where $\tau(t)>0 $ represents the time varying delay, $0<\tau_{0}\leq\tau(t)\leq\overline{\tau} $ and $\mu_1, $ $\mu_2 $ are positive constants. Under the assumptions $\frac{\rho_1}{K}=\frac{\rho_2}{b} $ and  $\mu_2<\sqrt{1-d}\mu_1, $ where $d$ is a constant such that $\tau'(t)\leq d<1, $ they proved that the energy decays exponentially.

In the presence of the thermo-viscoelastic damping, Djebabla and Tatar \cite{djebabla2010} considered the following Timoshenko system:
$$ \left\{
\begin{array}{lll}\displaystyle
\rho_1\varphi_{tt}-K(\varphi_{x}+\psi)_{x}=0, \ \ & (x,t)\in (0,L)\times(0,\infty), \medskip\\
\rho_2\psi_{tt}-b\psi_{xx}+K(\varphi_{x}+\psi)+\gamma\theta_{x}=0, \ \ & (x,t)\in (0,L)\times(0,\infty), \medskip\\
\displaystyle\rho_3\theta_{tt}-\delta\theta_{xx}+\beta\int_{0}^{t}g(t-s)\theta_{xx}(s){\rm d}s+\gamma\psi_{ttx}=0
, \ \ & (x,t)\in (0,L)\times(0,\infty),
\end{array}
\right.
$$
where $\rho_1, \rho_2, \rho_3, K, b, \beta, \gamma $ and $\delta $ are positive constants. They proved the exponential deacy of solutions in the energy norm if and only if the coefficients satisfy $\frac{b\rho_1}{K}-\rho_2=\delta-\frac{K\rho_3}{\rho_1}=\gamma $ and $g $ decays uniformly.

Kirane and Said-Houari \cite{kirane2011} examined the system of viscoelastic wave equations with a linear damping and a delay term
\bes
\left\{
\begin{array}{lll}\displaystyle
\displaystyle u_{tt}-\Delta u+\int_{0}^{t}g(t-s)\Delta u(x,s){\rm d}s+\mu_1u_{t}(x,t)\medskip\\
\quad\quad\qquad+\mu_2u_{t}\left(x,t-\tau\right)=0, \ \ &(x,t)\in \Omega\times(0,\infty), \medskip\\
u(x,t)=0, \ \ &(x,t)\in \partial\Omega\times(0,\infty), \medskip\\
u(x,0)=u_0(x), \ \ u_{t}(x,0)=u_1(x), \ \ &x\in \Omega, \medskip\\
u_{t}(x,t-\tau)=f_0(x,t-\tau), \ \ &(x,t)\in \Omega\times(0,\tau),
\end{array}
\right.
\lbl{6.2}
\ees
where $\Omega $ is a regular and bounded domain of $\mathbb{R}^{N} (N\geq1),  $ $\mu_{1}, \mu_2 $ are positive constants, $\tau>0 $ represents the time delay and $u_0, u_1, f_0 $ are given functions belonging to suitable spaces. They proved that the energy of problem \eqref{6.2} decreases exponentially as $t $ tends to infinity provided that $0<\mu_2\leq\mu_1 $ and $g $ decays exponentially. Liu \cite{liu2013} considered the viscoelastic wave equation with a linear damping and a time-varying delay term in the feedback. Under suitable assumptions, he established a general decay result.

Motivated by above research, we consider the well-posedness and the general energy decay for problem \eqref{1.1}. First, using the Faedo-Galerkin approximations together with some energy estimates, and under some restriction on the parameters $\mu_1 $ and $\mu_2, $ the system is showed to be well-posed.  Then, under the hypothesis $\mu_2\leq\mu_{1}, $ we prove a general decay of the total energy of our problem by using energy method. Our method of proof uses some ideas developed in \cite{kirane2011} for the wave equation with a viscoelastic damping and a delay term, enabling us to obtain suitable Lyapunov functionals, from which we derive the desired results. We recall that for $\mu_1=\mu_2, $ Nicaise and Pignotti showed in \cite{nicaise2006} that some instabilities may occur. Here, due to the presence of the viscoelastic damping, we prove that our energy still decays generally even if $\mu_1=\mu_2$.

The remaining part of this paper is organized as follows. In Section \ref{s2}, we present some materials and recall some useful lemmas needed for our work and state our main results. In Section \ref{s4}, we will prove the well-posedness of the solution. We will prove several technical lemmas and the general decay result under the two cases: $\mu_2<\mu_{1} $ and $\mu_2=\mu_{1} $ in Section \ref{s3}.

\section{Preliminaries and main results}\label{s2}
In this section, we present some assumptions and state the main results. We use the standard Lebesgue space $L^2(0,1)$ and the Sobolev space $H_0^1(0,1)$ with their usual scalar produces and norms and define the following space $X $ as
$$X=\left[H_{0}^{1}(0,1)\times L^2(0,1)\right]^2\times V \times H,$$
where
$$V=H^{1}(0,1)\cap H$$
and
$$H=\{\theta\in L^2(0,1)|\theta_x(0,t)=\theta_x(1,t)=0, \forall t\in [0,\infty)\}.$$

First, in order to exhibit the dissipative nature of system \eqref{1.1}, as in \cite{messaoudi2008}, we introduce the new variables $\Phi=\varphi_{t} $ and $\Psi=\psi_{t}. $ Then, as in \cite{nicaise2008}, we introduce the function
$$z(x,\rho,t)=\theta_{t}(x,t-\rho\tau), \ \ (x,\rho,t)\in (0,1)\times(0,1)\times(0,\infty).$$
Then, we have
\bess
\tau z_{t}(x,\rho,t)+z_{\rho}(x,\rho,t)=0,  \ \ (x,\rho,t)\in (0,1)\times(0,1)\times(0,\infty).
\eess
Therefore, problem \eqref{1.1} is equivalent to
\be \left\{
\begin{array}{lll}\displaystyle
\rho_1\Phi_{tt}-K(\Phi_{x}+\Psi)_{x}=0, \ \ &(x,t)\in (0,1)\times(0,\infty), \medskip\\
\rho_2\Psi_{tt}-b\Psi_{xx}+K(\Phi_{x}+\Psi)+\beta\theta_{tx}=0, \ \ &(x,t)\in (0,1)\times(0,\infty), \medskip\\
\displaystyle\rho_3\theta_{tt}-\delta\theta_{xx}+\gamma\Psi_{tx}+\int_{0}^{t}g(t-s)\theta_{xx}(s){\rm d}s\medskip\\
\quad\quad+\mu_{1}\theta_{t}(x,t)+\mu_{2}z(x,1,t)=0, \ \ &(x,t)\in (0,1)\times(0,\infty), \medskip\\
\tau z_{t}(x,\rho,t)+z_\rho(x,\rho,t)=0, \ \ &(x,\rho,t)\in (0,1)\times(0,1)\times(0,\infty)
\end{array}
\right.\lbl{1.5}
\ee
with the following initial datum and boundary conditions:
\be \left\{
\begin{array}{lll}\displaystyle
\Phi(x,0)=\Phi_0, \ \ \Phi_{t}(x,0)=\Phi_1, \ \ \Psi(x,0)=\Psi_0, \ \ \Psi_{t}(x,0)=\Psi_1, \ \ & x\in [0,1], \medskip\\
\theta(x,0)=\theta_0, \ \ \theta_{t}(x,0)=\theta_1, \ \ & x\in [0,1],  \medskip\\
\Phi(0,t)=\Phi(1,t)=\Psi(0,t)=\Psi(1,t)=\theta_{x}(0,t)=\theta_{x}(1,t)=0, \ \ & t\in[0,\infty), \medskip\\
z(x,0,t)=\theta_{t}(x,t), \ \ &(x,t)\in (0,1)\times(0,\infty), \medskip\\
\theta_{t}(x,t-\tau)=f_{0}(x,t-\tau), \ \ & (x,t)\in (0,1)\times(0,\tau).
\end{array}
\right.\lbl{12.5}
\ee

Next, we denote by $*$ the usual convolution term
$$(g*h)(t)=\int_{0}^{t}g(t-s)h(s){\rm d}s$$
and the binary operators $\diamond $ and $\circ $, respectively, by
$$(g\diamond h)(t)=\int_{0}^{t}g(t-s)\left(h(t)-h(s)\right){\rm d}s$$
and
$$(g\circ h)(t)=\int_{0}^{t}g(t-s)\left(h(t)-h(s)\right)^{2}{\rm d}s.$$

The following lemma was introduced in \cite{kirane2011}. It will be used in Section \ref{s3} to prove the general energy decay result for problem \eqref{1.1}-\eqref{1.11}.
\begin{lemma}\label{31}
For any function $g\in C^1(\mathbb{R})$ and any $h\in H^{1}(0,1),$ we have
\bes
\left(g*h\right)(t)h_t(t)&=&-\frac{1}{2}g(t)|h(t)|^2+\frac{1}{2}\left(g'\diamond h\right)(t)\nonumber\\
&&-\frac{1}{2}\frac{{\rm d}}{{\rm d}t}\left\{\left(g\diamond h\right)(t)-\left(\int_{0}^{t}g(s){\rm d}s\right)|h(t)|^2\right\}.
\lbl{31.4}
\ees
\end{lemma}
The proof of this lemma follows by differentiating the term $g\diamond h$.

\begin{lemma}\label{1}
(\cite{djebabla2010}) For any function $g\in C([0,\infty),\mathbb{R_{+}})$ and any $h\in L^{2}(0,1),$ we have
\bes
[(g\diamond h)(t)]^{2}\leq\left(\int_{0}^{t}g(s){\rm d}s\right)\left(g\circ h\right)(t), \ \ t\geq0.
\lbl{1.4}
\ees
\end{lemma}

Now, we assume that the kernel $g $ satisfies the following assumptions:

$(H1)$ $g: \mathbb{R_{+}}\rightarrow \mathbb{R_{+}}$ is a differential function such that
$$g(0)>0, \ \ \lambda=\delta-\int_{0}^{\infty}g(s){\rm d}s=\delta-\bar{g}>0;$$

$(H2)$ There exists a non-increasing differential function $\zeta: \mathbb{R_{+}}\rightarrow \mathbb{R_{+}}$ satisfying
$$g'(t)\leq -\zeta(t)g(t), \ \ t\geq0$$
and
$$\int_{0}^{+\infty}\zeta(t){\rm d}t=+\infty.$$

To state our decay result, we introduce the energy functional associated to problem \eqref{1.5}
\bes
E(t)&=&\frac{\gamma}{2}\int_{0}^{1}\left\{\rho_1\Phi_{t}^{2}+\rho_2\Psi_{t}^{2}+K|\Phi_{x}+\Psi|^{2}+b\Psi_{x}^{2}\right\}{\rm d}x\nonumber\\
&&+\frac{\beta}{2}\int_{0}^{1}\left\{\rho_3\theta_{t}^{2}+\left(\delta-\int_{0}^{t}g(s){\rm d}s\right)\theta_{x}^{2}+(g\circ\theta_{x})+\xi\int_{0}^{1}z^{2}(x,\rho,t){\rm d}\rho\right\}{\rm d}x,
\lbl{1.7}
\ees
where $\xi $ is a positive constant such that
\bes
\tau\mu_2<\xi<\tau(2\mu_1-\mu_2) \ \ {\rm if} \ \  \mu_2<\mu_1, \lbl{100}\medskip\\
\xi=\tau\mu_2 \ \ {\rm if} \ \  \mu_2=\mu_1.
\lbl{1.6}
\ees

Our main results read as follows.
\begin{theo}\label{32}
Assume that $\mu_2\leq\mu_1. $ Then for any given $\left(\Phi_{0}, \Phi_{1}, \Psi_{0}, \Psi_{1}, \theta_{0}, \theta_{1}\right)\in X,$ $\ f_0\in L^{2}\left((0,1)\times(0,1)\right)$ and $T>0, $ there exists a unique weak solution $\left(\Phi,\Psi,\theta,z\right)$ of problem \eqref{1.5}-\eqref{12.5} on $(0,T)$ such that
$$\left(\Phi, \Psi, \theta\right)\in C\left([0,T],\left[H^1_{0}(0,1)\right]^2\times V\right)\cap C^1\left([0,T],\left[L^2(0,1)\right]^2\times H\right).$$
\end{theo}

\begin{theo}\label{3}
Assume that $\mu_2\leq\mu_1 $ and $g$ satisfies $(H1) $ and $(H2). $ Assume further that initial datum satisfy
\bes
\left(\Phi_{0}, \Phi_{1}, \Psi_{0}, \Psi_{1}, \theta_{0}, \theta_{1}\right)\in X, \ f_0\in L^{2}\left((0,1)\times(0,1)\right)
\lbl{2.2}
\ees
and the coefficients $\rho_1, \rho_2, \rho_3, K, b, \gamma, \delta $ and $\beta$ satisfy
\bes
\frac{b\rho_1}{K}-\rho_2=\gamma,\ \ \delta-\frac{K\rho_3}{\rho_1}=\beta.
\lbl{2.3}
\ees
Then for any $t_{0}>0, $ there exist two positive constants $A $ and $\omega $ independent of the initial datum such that
\bes
E(t)\leq A e^{-\omega\int_{t_{0}}^{t}\zeta(s){\rm d}s}, \ \ t\geq t_{0}.
\lbl{2.4}
\ees
\end{theo}

\section{Proof of Theorem \ref{32} }\label{s4}
In this section, we will use the Faedo-Galerkin approximations together with some energy estimates, to prove the existence of the unique solution of problem \eqref{1.5}-\eqref{12.5} as stated in Theorem \ref{32}.

\noindent{\bf Proof.}  We divide the proof of Theorem \ref{32} in two steps: the construction of approximations and then thanks to certain energy estimates, we pass to the limit.

{\large\textbf{Step 1: Fadeo-Galerkin approximations.}}

As in \cite{kirane2011} and \cite{djebablat2011}, we construct approximations of the solution $(\Phi,\Psi,\theta,z)$ by the Faedo-Galerkin method as follows. For every $n\geq1,$ let $W_n$=span$\{\omega_1,\ldots,\omega_n\}$ be a Hilbet basis of the space $H_{0}^{1}(0,1).$

Now, we define for $1\leq j \leq n$ the sequence $\overline{\varphi}_{j}(x,\rho)$ as follows
$$\overline{\varphi}_{j}(x,0)=\omega_j(x).$$
Then, we may extend $\overline{\varphi}_{j}(x,0)$ by $\overline{\varphi}_{j}(x,\rho)$ over $L^2\left((0,1)\times[0,1]\right)$ and denote $V_n$=span$\{\overline{\varphi}_1,\ldots,\overline{\varphi}_n\}.$ We choose sequences $\left(\Phi_{0n},\Psi_{0n},\theta_{0n}\right)$ and $\left(\Phi_{1n},\Psi_{1n},\theta_{1n}\right)$ in $W_n$ and a sequence $\left(z_{0n}\right)$ in $V_n$ such that $\left(\Phi_{0n}, \Phi_{1n}, \Psi_{0n}, \Psi_{1n}, \theta_{0n}, \theta_{1n}\right)\rightarrow \left(\Phi_0, \Phi_1, \Psi_0, \Psi_1, \theta_0, \theta_1\right)$ strongly in $X$ and $z_{0n} \rightarrow f_0$ strongly in $L^2\left((0,1)\times(0,1)\right).$

We define now the approximations:
\bes
\left(\Phi_n(x,t),\Psi_n(x,t),\theta_n(x,t)\right)=\sum_{j=1}^{n}\left(f_{jn}(t),y_{jn}(t),h_{jn}(t)\right)\omega_j(x)
\ees
and
\bes
z_n(x,\rho,t)=\sum_{j=1}^{n}l_{jn}(t)\overline{\varphi}_{j}(x,\rho),
\ees
where $\left(\Phi_{n}(t),\Psi_{n}(t),\theta_{n}(t),z_{n}(t)\right) $ satisfies the following problem:
\be \left\{
\begin{array}{lll}\displaystyle
\rho_1\int_{0}^{1}\Phi_{ttn}\omega_j{\rm d}x+K\int_{0}^{1}(\Phi_{xn}+\Psi_n)\omega_{xj}{\rm d}x=0, \medskip\\
\displaystyle\rho_2\int_{0}^{1}\Psi_{ttn}\omega_j{\rm d}x+b\int_{0}^{1}\Psi_{xn}\omega_{xj}{\rm d}x+K\int_{0}^{1}(\Phi_{xn}+\Psi_{n})\omega_j{\rm d}x-\beta\int_{0}^{1}\theta_{tn}\omega_{xj}{\rm d}x=0, \medskip\\
\displaystyle\rho_3\int_{0}^{1}\theta_{ttn}\omega_j{\rm d}x+\delta\int_{0}^{1}\theta_{xn}\omega_{xj}{\rm d}x-\gamma\int_{0}^{1}\Psi_{tn}\omega_{xj}{\rm d}x-\int_{0}^{t}g(t-s)\int_{0}^{1}\theta_{xn}(s)\omega_{xj}{\rm d}x{\rm d}s\medskip\\
\quad\quad\quad\quad\quad\quad\quad\displaystyle+\int_{0}^{1}\left(\mu_{1}\theta_{tn}(x,t)+\mu_{2}z_n(x,1,t)\right)\omega_j{\rm d}x=0, \medskip\\
\left(\Phi_{n}(0), \Psi_{n}(0), \theta_{n}(0)\right)=\left(\Phi_{0n}, \Psi_{0n}, \theta_{0n} \right),\medskip\\
\left(\Phi_{tn}(0), \Psi_{tn}(0), \theta_{tn}(0)\right)=\left(\Phi_{1n}, \Psi_{1n}, \theta_{1n} \right),\medskip\\
z_{n}(x,0,t)=\theta_{tn}(x,t)
\end{array}
\right.\lbl{31.5}
\ee
and
\be \left\{
\begin{array}{lll}\displaystyle
\int_{0}^{1}\left(\tau z_{tn}(x,\rho,t)+z_{\rho n}(x,\rho,t)\right)\overline{\varphi}_{j}{\rm d}x=0,\medskip\\
z_n(x,\rho,0)=z_{0n},
\end{array}
\right.\lbl{31.6}
\ee
for $1\leq j \leq n.$ According to the standard theory of ordinary differential equations, the finite dimensional problem \eqref{31.5}-\eqref{31.6} has a solution $\left(f_{jn}(t),y_{jn}(t),h_{jn}(t),l_{jn}(t)\right)_{j=1,\ldots,n}$ defined on $[0,t_n).$ Then a priori estimates that follow imply that in fact $t_n=T.$

{\large\textbf{Step 2: Energy estimates.}}

Multiplying Eq. $\eqref{1.5}_{1}$ by $\gamma f'_{jn}, $ $\eqref{1.5}_{2}$ by $\gamma y'_{jn} $ and $\eqref{1.5}_{3}$ by $\beta  h'_{jn}, $ integrating over $(0,1) $ using integration by parts and Lemma \ref{31}, we get, for every $n\geq1,$
\bes
&&\frac{\gamma}{2}\left[\rho_1\|\Phi_{tn}\|_2^{2}+\rho_2\|\Psi_{tn}\|_2^{2}+K\|\Phi_{xn}+\Psi_n\|_2^{2}+b\|\Psi_{xn}\|_2^{2}\right]\nonumber\\
&&+\frac{\beta}{2}\left[\rho_3\|\theta_{tn}\|_2^{2}+\left(\delta-\int_{0}^{t}g(s){\rm d}s\right)\|\theta_{xn}\|_2^{2}+(g\circ\theta_{xn})\right]\nonumber\\
&&+\beta\mu_1\int_{0}^{1}\|\theta_{tn}(s)\|^2_2{\rm d}s+\beta\mu_2\int_{0}^{t}\int_{0}^{1}\theta_{tn}(x,s)z_n(x,1,s){\rm d}x{\rm d}s+\frac{\beta}{2}\int_{0}^{t}g(s)\|\theta_{xn}(s)\|_2^{2}{\rm d}s\nonumber\\
&&-\frac{\beta}{2}\int_{0}^{t}(g'\circ\theta_{xn})(s){\rm d}s\nonumber\\
&=&\frac{\gamma}{2}\left[\rho_1\|\Phi_1\|_2^2+\rho_2\|\Psi_1\|_2^2+K\|\Phi_{x0}+\Psi_0\|_2^2+b\|\Psi_{x0}\|_2^2\right]+\frac{\beta}{2}\left[\rho_3\|\theta_1\|_2^2+\delta\|\theta_{x0}\|_2^2\right].
\lbl{31.7}
\ees

Let $\xi>0$ to be chosen later. Multiplying Eq. $\eqref{1.5}_{4}$ by $\frac{\xi}{\tau}l'_{jn}(t) $ integrating over $(0,t)\times (0,1), $ we obtain
\bes
&&\frac{\xi}{2}\int_{0}^{1}\int_{0}^{1}z_{n}^2(x,\rho,t){\rm d}\rho{\rm d}x+\frac{\xi}{\tau}\int_{0}^{t}\int_{0}^{1}\int_{0}^{1}z_{n\rho}z_{n}(x,\rho,s){\rm d}\rho{\rm d}x{\rm d}s\nonumber\\
&=&\frac{\xi}{2}\|z_{0n}\|_{L^{2}\left((0,1)\times(0,1)\right)}^2.
\lbl{31.8}
\ees
Now, to handle the last term in the left-hand side of \eqref{31.8}, we remark that
\bes
\int_{0}^{t}\int_{0}^{1}\int_{0}^{1}z_{n\rho}z_{n}(x,\rho,t){\rm d}\rho{\rm d}x{\rm d}s&=&\frac{1}{2}\int_{0}^{t}\int_{0}^{1}\int_{0}^{1}\frac{\partial}{\partial \rho}z_{n}^2(x,\rho,s){\rm d}\rho{\rm d}x{\rm d}s\nonumber\\
&=&\frac{1}{2}\int_{0}^{t}\int_{0}^{1}\left(z_{n}^2(x,1,s)-z^{2}_n(x,0,s)\right){\rm d}x{\rm d}s.
\lbl{31.9}
\ees
Summing up the identities \eqref{31.7} and \eqref{31.8} and taking into account \eqref{31.9}, we get
\bes
&&\mathscr{E}_n(t)+\beta\left(\mu_1-\frac{\xi}{2\tau}\right)\int_{0}^{t}\|\theta_{tn}\|^2_2{\rm d}s+\frac{\xi\beta}{2\tau}\int_{0}^{t}\int_{0}^{1}z_{n}^2(x,1,s){\rm d}x{\rm d}s\nonumber\\
&&+\beta\mu_2\int_{0}^{t}\int_{0}^{1}z_n(x,1,s)\theta_{tn}(x,s){\rm d}s\nonumber\\
&&+\frac{\beta}{2}\int_{0}^{t}g(s)\|\theta_{xn}(s)\|^2_2{\rm d}s-\frac{\beta}{2}\int_{0}^{t}\left(g'\circ\theta_{xn}\right)(s){\rm d}s\nonumber\\
&=&\mathscr{E}_n(0),
\lbl{40.1}
\ees
where
\bes
\mathscr{E}_n(t)&=&\frac{\gamma}{2}\left[\rho_1\|\Phi_{tn}\|_2^{2}+\rho_2\|\Psi_{tn}\|_2^{2}+K\|\Phi_{xn}+\Psi_{n}\|_2^{2}+b\|\Psi_{xn}\|_2^{2}\right]\nonumber\\
&&+\frac{\beta}{2}\left[\rho_3\|\theta_{tn}\|_2^{2}+\left(\delta-\int_{0}^{t}g(s){\rm d}s\right)\|\theta_{xn}\|_2^{2}+(g\circ\theta_{xn})+\xi\|z_n\|_{L^{2}{\left((0,1)\times(0,1)\right)}}^{2}\right].
\lbl{40.0}
\ees
At this point, we have to distinguish the following two cases:

Case 1: We suppose that $\mu_2<\mu_1. $ Let us choose $\xi $ satisfies inequality \eqref{100}. Using Young's inequality, \eqref{40.1} leads to:
\bess
&&\mathscr{E}_n(t)+\beta\left(\mu_1-\frac{\xi}{2\tau}-\frac{\mu_2}{2}\right)\int_{0}^{t}\|\theta_{tn}\|^2_2{\rm d}s+\beta\left(\frac{\xi}{2\tau}-\frac{\mu_2}{2}\right)\int_{0}^{t}\int_{0}^{1}z_{n}^2(x,1,s){\rm d}x{\rm d}s\nonumber\\
&&+\frac{\beta}{2}\int_{0}^{t}g(s)\|\theta_{xn}(s)\|^2_2{\rm d}s-\frac{\beta}{2}\int_{0}^{t}\left(g'\circ\theta_{xn}\right)(s){\rm d}s\nonumber\\
&\leq&\mathscr{E}_n(0).
\eess
Consequently, using \eqref{1.6}, we can find two positive constants $c_1$ and $c_2$ such that
\bes
&&\mathscr{E}_n(t)+c_1\int_{0}^{t}\|\theta_{tn}\|^2_2{\rm d}s+c_2\int_{0}^{t}\int_{0}^{1}z_{n}^2(x,1,s){\rm d}x{\rm d}s\nonumber\\
&&+\frac{\beta}{2}\int_{0}^{t}g(s)\|\theta_{xn}(s)\|^2_2{\rm d}s-\frac{\beta}{2}\int_{0}^{t}\left(g'\circ\theta_{xn}\right)(s){\rm d}s\nonumber\\
&\leq&\mathscr{E}_n(0).
\lbl{40.2}
\ees

Case 2: We suppose that $\mu_2=\mu_1=\mu$ and choose $\xi=\tau\mu.$ Then, inequality \eqref{40.2} takes the form
\bes
\mathscr{E}_n(t)+\frac{\beta}{2}\int_{0}^{t}g(s)\|\theta_{xn}(s)\|^2_2{\rm d}s-\frac{\beta}{2}\int_{0}^{t}\left(g'\circ\theta_{xn}\right)(s){\rm d}s\leq\mathscr{E}_n(0).
\lbl{40.3}
\ees
Now, in both cases and since the sequences $(\Phi_{0n})_{n\in N},$ $(\Phi_{1n})_{n\in N},$ $(\Psi_{0n})_{n\in N},$ $(\Psi_{1n})_{n\in N},$ $(\theta_{0n})_{n\in N},$ $(\theta_{1n})_{n\in N}$ and $(z_{0n})_{n\in N} $ converge, and using $(H1)$ and $(H2),$ we can find a positive constant $C$ independent of $n$ such that
\bes
\mathscr{E}_n(t)\leq C.
\lbl{40.4}
\ees
Therefore, using the fact that $\delta-\int_{0}^{1}g(s){\rm d}s\geq\lambda, $ the last estimate \eqref{40.4} together with \eqref{40.0} give us, for all $n\in\mathbb{N}, t_{n}=T, $ we deduce that
\bes
(\Phi_n,\Psi_n,\theta_n)_{n\in\mathbb{N}}\ \  {\rm is}\  \  {\rm bounded} \ \ {\rm in} \ \
L^{\infty}\left(0,T;\left[H^1_{0}(0,1)\right]^2\times V\right),
\lbl{40.5}
\ees
\bes
(\Phi_{tn},\Psi_{tn},\theta_{tn})_{n\in\mathbb{N}} \ \  {\rm is}\  \  {\rm bounded} \ \ {\rm in} \ \
L^{\infty}\left(0,T;\left[L^2(0,1)\right]^2\times H\right)
\lbl{40.6}
\ees
and
\bes
(z_n)_{n\in\mathbb{N}} \ \  {\rm is}\  \  {\rm bounded} \ \ {\rm in} \ \
L^{\infty}\left(0,T;L^{2}\left((0,1)\times(0,1)\right)\right).
\lbl{40.7}
\ees
Consequently, we may conclude that
\bes
\left(\Phi_n,\Psi_n,\theta_n\right)\rightharpoonup\left(\Phi,\Psi,\theta\right) \   {\rm weak^{*}}\ \ \quad {\rm in} \ \ L^{\infty}\left(0,T;\left[H^1_{0}(0,1)\right]^2\times V\right),\nonumber\\
\left(\Phi_{tn},\Psi_{tn},\theta_{tn}\right)\rightharpoonup\left(\Phi_t,\Psi_t,\theta_t\right) \   {\rm weak^{*}}\ \ \quad {\rm in} \ \
 L^{\infty}\left(0,T;\left[L^2(0,1)\right]^2\times H\right),\nonumber\\
z_n\rightharpoonup z \   {\rm weak^{*}}\ \ \quad {\rm in} \ \  L^{\infty}\left(0,T;L^{2}((0,1)\times(0,1))\right).
\ees
From \eqref{40.5}, \eqref{40.6} and \eqref{40.7}, we have $(\Phi_n,\Psi_n,\theta_n)_{n\in\mathbb{N}}$ is bounded in $L^{\infty}\left(0,T;\left[H^1_{0}(0,1)\right]^2\times V\right).$ Then, $(\Phi_n,\Psi_n,\theta_n)_{n\in\mathbb{N}}$ is bounded in $L^{2}\left(0,T;\left[H^1_{0}(0,1)\right]^2\times V\right).$  Since $(\Phi_{tn},\Psi_{tn},\theta_{tn})_{n\in\mathbb{N}}$ is bounded in $L^{\infty}\left(0,T;\left[L^2(0,1)\right]^2\times H\right),$ $(\Phi_{tn},\Psi_{tn},\theta_{tn})_{n\in\mathbb{N}} $ is bounded in $L^{2}\left(0,T;\left[L^2(0,1)\right]^2\times H\right).$ Consequently, $(\Phi_n,\Psi_n,\theta_n)_{n\in\mathbb{N}}$ is bounded in $H^{1}\left(0,T;\left[H^1(0,1)\right]^2\times V\right).$

Since the embedding $H^{1}\left(0,T;\left[H^1(0,1)\right]^2\times V\right)\hookrightarrow L^{2}\left(0,T;\left[L^2(0,1)\right]^2\times H\right)$ is compact, using Aubin-Lions theorem \cite{l1969}, we can extract subsequences $(\Phi_\mu,\Psi_\mu,\theta_\mu)_{\mu\in\mathbb{N}}$ of $(\Phi_n,\Psi_n,\theta_n)_{n\in\mathbb{N}}$   such that
\bess
\left(\Phi_\mu,\Psi_\mu,\theta_\mu\right)\rightharpoonup\left(\Phi,\Psi,\theta\right) \   {\rm strongly}\  {\rm in} \ \ L^{2}\left(0,T;\left[L^2(0,1)\right]^2\times H\right).
\eess
Therefore,
\bess
\left(\Phi_\mu,\Psi_\mu,\theta_\mu\right)\rightharpoonup\left(\Phi,\Psi,\theta\right) \   {\rm strongly}\ \ {\rm and} \  \ {\rm a.e}\  \ {\rm on} \ \ (0,T)\times(0,1).
\eess
The proof now can be completed arguing as in \cite{l1969}.

\section{Proof of Theorem \ref{3}}\label{s3}
In this section, under the hypothesis $\mu_2\leq\mu_1, $ we show that the energy of the solution of problem \eqref{1.5}-\eqref{12.5} decreases generally as $t $ tends to infinity by using the energy method and suitable Lyapunov functionals. We will separately discuss the two cases which are the case $\mu_2<\mu_1 $ and the case $\mu_2=\mu_1 $ since the proofs are slightly different.

{\large\textbf{4.1 The case $\mu_2<\mu_1$}}

Our goal now is to prove the above energy $E(t) $ is a non-increasing functional along the trajectories. More precisely, we have the following result:
\begin{lemma}\lbl{2}
Suppose that $(H1) $ and $(H2) $ hold and let $(\Phi,\Psi,\theta,z)$ be a solution of \eqref{1.5}-\eqref{12.5}. Then we have
\bes
E'(t)&\leq&-\frac{\beta}{2}g(t)\int_{0}^{1}\theta_{x}^{2}{\rm d}x+\frac{\beta}{2}\int_{0}^{1}(g'\circ\theta_{x}){\rm d}x-\beta\left(\mu_1-\frac{\xi}{2\tau}-\frac{\mu_2}{2}\right)\int_{0}^{1}\theta_{t}^{2}{\rm d}x\nonumber\\
&&-\beta\left(\frac{\xi}{2\tau}-\frac{\mu_2}{2}\right)\int_{0}^{1}z^{2}(x,1,t){\rm d}x\leq 0, \ \ t\geq0.
\lbl{1.8}
\ees
\end{lemma}
{\bf Proof.\  } Multiplying Eq. $\eqref{1.5}_{1}$ by $\gamma\Phi_{t}, $ $\eqref{1.5}_{2}$ by $\gamma\Psi_{t} $ and $\eqref{1.5}_{3}$ by $\beta\theta_{t}, $ integrating over $(0,1) $ and performing an integration by parts, we find
\bes
&&\frac{\gamma}{2}\frac{{\rm d}}{{\rm d}t}\left[\int_{0}^{1}\left\{\rho_1\Phi_{t}^{2}+\rho_2\Psi_{t}^{2}+K|\Phi_{x}+\Psi|^{2}+b\Psi_{x}^{2}\right\}{\rm d}x\right]\nonumber\\
&&+\frac{\beta}{2}\frac{{\rm d}}{{\rm d}t}\left[\int_{0}^{1}\left\{\rho_3\theta_{t}^{2}+\left(\delta-\int_{0}^{t}g(\tau){\rm d}\tau\right)\theta_{x}^{2}+(g\circ\theta_{x})\right\}{\rm d}x\right]\nonumber\\
&=&-\beta\mu_1\int_{0}^{1}\theta^{2}_{t}{\rm d}x-\beta\mu_2\int_{0}^{1}\theta_{t}z(x,1,t){\rm d}x-\frac{\beta}{2}g(t)\int_{0}^{1}\theta_{x}^{2}{\rm d}x+\frac{\beta}{2}\int_{0}^{1}(g'\circ\theta_{x}){\rm d}x.
\lbl{1.9}
\ees
Now, multiplying the Eq. $\eqref{1.5}_{4}$ by $\frac{\xi}{\tau}z, $ integrating the result over $(0,1)\times (0,1) $ with respect to $\rho $ and $x,$ respectively, we obtain
\bes
\frac{\xi}{2}\frac{{\rm d}}{{\rm d}t}\int_{0}^{1}\int_{0}^{1}z^{2}(x,\rho,t){\rm d}\rho{\rm d}x&=&-\frac{\xi}{\tau}\int_{0}^{1}\int_{0}^{1}z z_{\rho}(x,\rho,t){\rm d}\rho{\rm d}x\nonumber\\
&=&-\frac{\xi}{2\tau}\int_{0}^{1}\int_{0}^{1}\frac{\partial}{\partial\rho}z^{2}(x,\rho,t){\rm d}\rho{\rm d}x\nonumber\\
&=&\frac{\xi}{2\tau}\int_{0}^{1}\left[z^{2}(x,0,t)-z^{2}(x,1,t)\right]{\rm d}x\nonumber\\
&=&\frac{\xi}{2\tau}\int_{0}^{1}\left(\theta_{t}^{2}-z^{2}(x,1,t)\right){\rm d}x.
\lbl{2.0}
\ees
By Young's inequality, we have
\bes
-\mu_2\int_{0}^{1}\theta_{t}z(x,1,t){\rm d}x\leq\frac{\mu_2}{2}\int_{0}^{1}\theta_{t}^{2}{\rm d}x+\frac{\mu_2}{2}\int_{0}^{1}z^{2}(x,1,t){\rm d}x.
\lbl{2.1}
\ees
Then, exploiting \eqref{2.0} and \eqref{2.1} our conclusion holds. \quad \quad $\Box$

Now we are going to construct a Lyapunov functional $\mathcal{L}$ equivalent to $E.$ For this, we define several functionals which allow us to obtain the needed estimates. We introduce the multiplier $w$ given by the solution of the Dirichlet problem
\bes
-w_{xx}=\Psi_{x},\ \ w(0)=w(1)=0
\lbl{2.5}
\ees
and  define the functional
\bes
I_{1}(t)=\int_{0}^{1}\left(\rho_2\Psi_{t}\Psi+\rho_1\Phi_{t}w+\beta\theta_{x}\Psi\right){\rm d}x,\ \ t\geq0.
\lbl{2.6}
\ees
The derivative of this functional will provide us with the term
$$-\int_{0}^{1}\Psi_{x}^{2}{\rm d}x.$$

\begin{lemma}\label{4}
Let $(\Phi,\Psi,\theta,z)$ be a solution of \eqref{1.5}-\eqref{12.5}. For any $\varepsilon_{1}>0, $ we have
\bes
I'_{1}(t)\leq-b\int_{0}^{1}\Psi_{x}^{2}{\rm d}x+\left(\frac{3\rho_2}{2}+\frac{\rho_1^{2}C_{p}}{4\varepsilon_{1}}\right)\int_{0}^{1}\Psi_{t}^{2}{\rm d}x+\varepsilon_{1}\int_{0}^{1}\Phi_{t}^{2}{\rm d}x+\frac{\beta^{2}}{2\rho_2}\int_{0}^{1}\theta_{x}^{2}{\rm d}x,
\lbl{2.7}
\ees
where $C_{p}$ is the Poincar$\acute{e}$ constant.
\end{lemma}
{\bf Proof. } By differentiating $I_1$ with respect to $t$ and using Eqs. $\eqref{1.5}$ and  \eqref{12.5} we conclude that
$$I'_{1}(t)=-b\int_{0}^{1}\Psi_{x}^{2}{\rm d}x+\rho_2\int_{0}^{1}\Psi_{t}^{2}{\rm d}x-K\int_{0}^{1}\Psi^{2}{\rm d}x+K\int_{0}^{1}w_{x}^{2}{\rm d}x+\rho_1\int_{0}^{1}\Phi_{t}w_{t}{\rm d}x+\beta\int_{0}^{1}\theta_{x}\Psi_{t}{\rm d}x.$$
By exploiting the inequalities
$$\int_{0}^{1}w_{x}^{2}{\rm d}x\leq\int_{0}^{1}\Psi^{2}{\rm d}x\leq C_{p}\int_{0}^{1}\Psi_{x}^{2}{\rm d}x,$$

$$\int_{0}^{1}w_{t}^{2}{\rm d}x\leq C_{p}\int_{0}^{1}w_{tx}^{2}{\rm d}x\leq C_{p}\int_{0}^{1}\Psi_{t}^{2}{\rm d}x$$
and Young's inequality, we find that
$$I'_{1}(t)\leq-b\int_{0}^{1}\Psi_{x}^{2}{\rm d}x+\frac{3\rho_2}{2}\int_{0}^{1}\Psi_{t}^{2}{\rm d}x+\frac{\beta^{2}}{2\rho_2}\int_{0}^{1}\theta_{x}^{2}{\rm d}x+\varepsilon_{1}\int_{0}^{1}\Phi_{t}^{2}{\rm d}x+\frac{\rho_1^{2}}{4\varepsilon_1}\int_{0}^{1}w_{t}^{2}{\rm d}x.$$
Thus,
$$I'_{1}(t)\leq-b\int_{0}^{1}\Psi_{x}^{2}{\rm d}x+\left(\frac{3\rho_2}{2}+\frac{\rho_1^{2}C_{p}}{4\varepsilon_{1}}\right)\int_{0}^{1}\Psi_{t}^{2}{\rm d}x+\varepsilon_{1}\int_{0}^{1}\Phi_{t}^{2}{\rm d}x+\frac{\beta^{2}}{2\rho_2}\int_{0}^{1}\theta_{x}^{2}{\rm d}x,$$
which is exactly \eqref{2.7}.\quad \quad $\Box$

In order to obtain the negative term of $\int_{0}^{1}\theta_{x}^{2}{\rm d}x, $ we define the functional
$$I_{2}(t)=\int_{0}^{1}\left(\rho_3\theta_{t}\theta+\gamma\Psi_{x}\theta+\frac{\mu_1}{2}\theta^{2}\right){\rm d}x.$$
\begin{lemma}\label{5}
Let $(\Phi,\Psi,\theta,z)$ be a solution of \eqref{1.5}-\eqref{12.5}. For any $\varepsilon_{2}>0, $ we have
\bes
I'_{2}(t)&\leq&-\frac{\lambda}{2}\int_{0}^{1}\theta_{x}^{2}{\rm d}x+\left(\rho_3+\frac{\gamma^{2}}{4\varepsilon_{2}}\right)\int_{0}^{1}\theta_{t}^{2}{\rm d}x+\varepsilon_{2}\int_{0}^{1}\Psi_{x}^{2}{\rm d}x\nonumber\\
&&+\frac{1}{\lambda}\int_{0}^{t}g(s){\rm d}s\int_{0}^{1}\left(g\circ\theta_{x}\right){\rm d}x+\frac{\mu_{2}^{2}C_p}{\lambda}\int_{0}^{1}z^{2}(x,1,t){\rm d}x.
\lbl{2.8}
\ees
\end{lemma}
{\bf Proof. } A simple differentiation leads to
$$I'_{2}(t)=\rho_3\int_{0}^{1}\theta_{t}^{2}{\rm d}x+\rho_3\int_{0}^{1}\theta_{tt}\theta{\rm d}x+\gamma\int_{0}^{1}\Psi_{tx}\theta{\rm d}x+\gamma\int_{0}^{1}\Psi_{x}\theta_{t}{\rm d}x+\mu_1\int_{0}^{1}\theta_{t}\theta{\rm d}x.$$
By using Eq. $\eqref{1.5}_{3} $ and \eqref{12.5}, we arrive at
$$I'_{2}(t)=-\left(\delta-\int_{0}^{t}g(s){\rm d}s\right)\int_{0}^{1}\theta_{x}^{2}{\rm d}x+\rho_3\int_{0}^{1}\theta_{t}^{2}{\rm d}x+\gamma\int_{0}^{1}\Psi_{x}\theta_{t}{\rm d}x-\int_{0}^{1}\left(g\diamond\theta_{x}\right)\theta_{x}{\rm d}x-\mu_2\int_{0}^{1}z(x,1,t)\theta{\rm d}x, $$
the last three terms can be estimated, using Young's inequality and Lemma \ref{1}, as follows, for $\alpha>0,$
$$\gamma\int_{0}^{1}\Psi_{x}\theta_{t}{\rm d}x\leq\varepsilon_{2}\int_{0}^{1}\Psi_{x}^{2}{\rm d}x+\frac{\gamma^{2}}{4\varepsilon_2}\int_{0}^{1}\theta_{t}^{2}{\rm d}x,$$
$$-\mu_2\int_{0}^{1}z(x,1,t)\theta{\rm d}x\leq\alpha\int_{0}^{1}\theta_{x}^{2}{\rm d}x+\frac{\mu_{2}^{2}C_p}{4\alpha}\int_{0}^{1}z^{2}(x,1,t){\rm d}x$$
and
$$-\int_{0}^{1}\left(g\diamond\theta_{x}\right)\theta_{x}{\rm d}x\leq\alpha\int_{0}^{1}\theta_{x}^{2}{\rm d}x+\frac{1}{4\alpha}\int_{0}^{t}g(s){\rm d}s\int_{0}^{1}\left(g\circ\theta_{x}\right){\rm d}x.$$
We deduce that
\bes
I'_{2}(t)&\leq&-\left(\delta-\int_{0}^{t}g(s){\rm d}s-2\alpha\right)\int_{0}^{1}\theta_{x}^{2}{\rm d}x+\left(\rho_3+\frac{\gamma^{2}}{4\varepsilon_{2}}\right)\int_{0}^{1}\theta_{t}^{2}{\rm d}x\nonumber\\
&&+\varepsilon_{2}\int_{0}^{1}\Psi_{x}^{2}{\rm d}x
+\frac{1}{4\alpha}\int_{0}^{t}g(s){\rm d}s\int_{0}^{1}\left(g\circ\theta_{x}\right){\rm d}x+\frac{\mu_{2}^{2}C_p}{4\alpha}\int_{0}^{1}z^{2}(x,1,t){\rm d}x.
\ees
The choice of $\alpha=\frac{\lambda}{4}$ gives the result.\quad \quad $\Box$

In order to get the negative terms of $\int_{0}^{1}\Phi_{t}^{2}{\rm d}x $ and $\int_{0}^{1}\Psi_{t}^{2}{\rm d}x, $ we define the functional $$I_{3}(t)=-\rho_1\int_{0}^{1}\Phi_{t}\Phi{\rm d}x-\rho_2\int_{0}^{1}\Psi_{t}\Psi{\rm d}x.$$
\begin{lemma}\label{6}
Let $(\Phi,\Psi,\theta,z)$ be a solution of \eqref{1.5}-\eqref{12.5}. Then we have
\bes
I'_{3}(t)&\leq&-\rho_1\int_{0}^{1}\Phi_{t}^{2}{\rm d}x-\rho_2\int_{0}^{1}\Psi_{t}^{2}{\rm d}x+\frac{3b}{2}\int_{0}^{1}\Psi_{x}^{2}{\rm d}x\nonumber\\
&&+K\int_{0}^{1}|\Phi_{x}+\Psi|^{2}{\rm d}x+\frac{\beta^{2}}{2b}\int_{0}^{1}\theta_{t}^{2}{\rm d}x.
\lbl{2.9}
\ees
\end{lemma}
{\bf Proof. } A differentiation of $I_{3},$ taking into account \eqref{1.5}-\eqref{12.5}, gives
$$I'_{3}(t)=-\rho_1\int_{0}^{1}\Phi_{t}^{2}{\rm d}x-\rho_2\int_{0}^{1}\Psi_{t}^{2}{\rm d}x+b\int_{0}^{1}\Psi_{x}^{2}{\rm d}x+K\int_{0}^{1}|\Phi_{x}+\Psi|^{2}{\rm d}x-\beta\int_{0}^{1}\theta_{t}\Psi_{x}{\rm d}x.$$
Using Young's and Poincar$\acute{e}$'s inequalities for the last term, we  obtain \eqref{2.9}.\quad \quad $\Box$

In order to get the negative term of $\int_{0}^{1}|\Phi_{x}+\Psi|^{2}{\rm d}x, $ we define the functional
\bes
I_{4}(t)&=&\rho_{2}\int_{0}^{1}\Psi_{t}(\Phi_{x}+\Psi){\rm d}x+(\rho_2+\gamma)\int_{0}^{1}\Psi_{x}\Phi_{t}{\rm d}x+\rho_3\int_{0}^{1}\theta_{t}\Phi_{t}{\rm d}x\nonumber\\
&&+\left(\frac{K\rho_3}{\rho_1}+\beta\right)\int_{0}^{1}\theta_{x}\Phi_{x}{\rm d}x-\int_{0}^{1}\left(g*\theta_{x}\right)\Phi_{x}{\rm d}x.
\ees
\begin{lemma}\label{7}
Let $(\Phi,\Psi,\theta,z) $ be a solution of \eqref{1.5}-\eqref{12.5} and assume that \eqref{2.3} holds. Then, for any $\varepsilon_4>0, $ we get
\bes
I'_{4}(t)&\leq&\varepsilon_4\left(1+\frac{\rho_3K^{2}}{\rho_1^{2}b}\right)[\Phi_{x}^{2}(1)+\Phi_{x}^{2}(0)]+\frac{b^{2}}{4\varepsilon_4}[\Psi_{x}^{2}(1)
+\Psi_{x}^{2}(0)]+\frac{b\rho_3}{4\varepsilon_4}[\theta_{t}^{2}(1)+\theta_{t}^{2}(0)]\nonumber\\
&&-K\int_{0}^{1}|\Phi_{x}+\Psi|^{2}{\rm d}x+\varepsilon_4\int_{0}^{1}\Psi_{x}^{2}{\rm d}x+\frac{\delta^{2}+\mu_1^{2}}{4\varepsilon_4}\int_{0}^{1}\theta_{t}^{2}{\rm d}x+2\varepsilon_4\int_{0}^{1}\Phi_{t}^{2}{\rm d}x\nonumber\\
&&+\rho_{2}\int_{0}^{1}\Psi_{t}^{2}{\rm d}x+\varepsilon_4\int_{0}^{1}\Phi_{x}^{2}{\rm d}x+\frac{g^{2}(0)}{2\varepsilon_4}\int_{0}^{1}\theta_{x}^{2}{\rm d}x-\frac{g(0)}{2\varepsilon_4}\int_{0}^{1}\left(g'\circ \theta_{x}\right){\rm d}x\nonumber\\
&&+\frac{\mu_{2}^{2}}{4\varepsilon_4}\int_{0}^{1}z^{2}(x,1,t){\rm d}x.
\lbl{3.0}
\ees
\end{lemma}
{\bf Proof. } By differentiating the functional $I_{4}, $ using Eq. $\eqref{1.5}_{1}, \eqref{1.5}_{2} $ and $\eqref{1.5}_{3}$, we obtain
\bes
I'_{4}(t)&=&\int_{0}^{1}\left[b\Psi_{xx}-K\left(\Phi_{x}+\Psi\right)-\beta\theta_{tx}\right]\left(\Phi_{x}+\Psi\right){\rm d}x+\rho_2\int_{0}^{1}\Psi_{t}\left(\Phi_{x}+\Psi\right)_{t}{\rm d}x\nonumber\\
&&+(\rho_2+\gamma)\int_{0}^{1}\Psi_{tx}\Phi_{t}{\rm d}x+(\rho_2+\gamma)\int_{0}^{1}\Psi_{x}\Phi_{tt}{\rm d}x+\frac{K\rho_3}{\rho_1}\int_{0}^{1}\theta_{t}\left(\Phi_{x}+\Psi\right)_{x}{\rm d}x\nonumber\\
&&+\int_{0}^{1}\left(\delta\theta_{xx}-\gamma\Psi_{tx}-\mu_1\theta_{t}-\mu_2z(x,1,t)-g*\theta_{xx}\right)\Phi_{t}{\rm d}x+\left(\frac{K\rho_3}{\rho_1}+\beta\right)\int_{0}^{1}\theta_{tx}\Phi_{x}{\rm d}x\nonumber\\
&&+\left(\frac{K\rho_3}{\rho_1}+\beta\right)\int_{0}^{1}\theta_{x}\Phi_{tx}{\rm d}x-\int_{0}^{1}\left(g*\theta_{x}\right)_{t}\Phi_{x}{\rm d}x-\int_{0}^{1}\left(g*\theta_{x}\right)\Phi_{tx}{\rm d}x\nonumber\\
&=&b[\Phi_{x}\Psi_{x}]_{x=0}^{x=1}-\frac{b\rho_1}{K}\int_{0}^{1}\Psi_{x}\Phi_{tt}{\rm d}x-K\int_{0}^{1}|\Phi_{x}+\Psi|^{2}{\rm d}x-\beta\int_{0}^{1}\theta_{tx}\Phi_{x}{\rm d}x+\beta\int_{0}^{1}\theta_{t}\Psi_{x}{\rm d}x\nonumber\\
&&+\rho_2\int_{0}^{1}\Psi_{t}\Phi_{tx}{\rm d}x+\rho_2\int_{0}^{1}\Psi_{t}^{2}{\rm d}x+(\rho_2+\gamma)\int_{0}^{1}\Psi_{tx}\Phi_{t}{\rm d}x+(\rho_2+\gamma)\int_{0}^{1}\Psi_{x}\Phi_{tt}{\rm d}x\nonumber\\
&&+\frac{K\rho_3}{\rho_1}\int_{0}^{1}\theta_{t}\Psi_{x}{\rm d}x+\frac{K\rho_3}{\rho_1}[\Phi_{x}\theta_{t}]_{x=0}^{x=1}-\frac{K\rho_3}{\rho_1}\int_{0}^{1}\theta_{tx}\Phi_{x}{\rm d}x-\delta\int_{0}^{1}\theta_{x}\Phi_{tx}{\rm d}x-\gamma\int_{0}^{1}\Psi_{tx}\Phi_{t}{\rm d}x\nonumber\\
&&-\int_{0}^{1}\left(\mu_{1}\theta_{t}+\mu_{2}z(x,1,t)\right)\Phi_{t}{\rm d}x+\int_{0}^{1}(g*\theta_{x})\Phi_{tx}{\rm d}x+\left(\frac{K\rho_3}{\rho_1}+\beta\right)\int_{0}^{1}\theta_{tx}\Phi_{x}{\rm d}x\nonumber\\
&&+\left(\frac{K\rho_3}{\rho_1}+\beta\right)\int_{0}^{1}\theta_{x}\Phi_{tx}{\rm d}x-\int_{0}^{1}\left(g*\theta_{x}\right)_{t}\Phi_{x}{\rm d}x-\int_{0}^{1}\left(g*\theta_{x}\right)\Phi_{tx}{\rm d}x\nonumber\\
&=&b[\Phi_{x}\Psi_{x}]_{x=0}^{x=1}+\frac{K\rho_3}{\rho_1}[\Phi_{x}\theta_{t}]_{x=0}^{x=1}-K\int_{0}^{1}|\Phi_{x}+\Psi|^{2}{\rm d}x+\left(\beta+\frac{K\rho_3}{\rho_1}\right)\int_{0}^{1}\theta_{t}\Psi_{x}{\rm d}x+\rho_2\int_{0}^{1}\Psi_{t}^{2}{\rm d}x\nonumber\\
&&-\int_{0}^{1}\left(\mu_{1}\theta_{t}+\mu_{2}z(x,1,t)\right)\Phi_{t}{\rm d}x-\int_{0}^{1}\left(g*\theta_{x}\right)_{t}\Phi_{x}{\rm d}x+\left(\rho_2+\gamma-\frac{b\rho_1}{K}\right)\int_{0}^{1}\Psi_{x}\Phi_{tt}{\rm d}x\nonumber\\
&&+\left(\frac{K\rho_3}{\rho_1}+\beta-\delta\right)\int_{0}^{1}\theta_{x}\Phi_{tx}{\rm d}x.
\ees
By using \eqref{2.3}, we have
\bes
I'_{4}(t)&=&b[\Phi_{x}\Psi_{x}]_{x=0}^{x=1}+\frac{K\rho_3}{\rho_1}[\Phi_{x}\theta_{t}]_{x=0}^{x=1}-K\int_{0}^{1}|\Phi_{x}+\Psi|^{2}{\rm d}x+\delta\int_{0}^{1}\theta_{t}\Psi_{x}{\rm d}x\nonumber\\
&&+\rho_2\int_{0}^{1}\Psi_{t}^{2}{\rm d}x-\int_{0}^{1}\left(\mu_{1}\theta_{t}+\mu_{2}z(x,1,t)\right)\Phi_{t}{\rm d}x-\int_{0}^{1}\left(g*\theta_{x}\right)_{t}\Phi_{x}{\rm d}x.
\lbl{3.1}
\ees
Now we estimate the terms in the right side of \eqref{3.1}. Applying Young's and Poincar$\acute{e}$'s inequalities and Lemma \ref{1}, we obtain that for any $\varepsilon_4>0, $
$$\delta\int_{0}^{1}\theta_{t}\Psi_{x}{\rm d}x\leq\varepsilon_4\int_{0}^{1}\Psi_{x}^{2}{\rm d}x+\frac{\delta^{2}}{4\varepsilon_4}\int_{0}^{1}\theta_{t}^{2}{\rm d}x,$$
$$-\mu_{1}\int_{0}^{1}\theta_{t}\Phi_{t}{\rm d}x\leq\varepsilon_4\int_{0}^{1}\Phi_{t}^{2}{\rm d}x+\frac{\mu_1^{2}}{4\varepsilon_4}\int_{0}^{1}\theta_{t}^{2}{\rm d}x,$$
$$-\mu_{2}\int_{0}^{1}z(x,1,t)\Phi_{t}{\rm d}x\leq\varepsilon_4\int_{0}^{1}\Phi_{t}^{2}{\rm d}x+\frac{\mu_2^{2}}{4\varepsilon_4}\int_{0}^{1}z^{2}(x,1,t){\rm d}x,$$
$$b[\Phi_{x}\Psi_{x}]_{x=0}^{x=1}\leq\varepsilon_4[\Phi_{x}^{2}(1)+\Phi_{x}^{2}(0)]+\frac{b^{2}}{4\varepsilon_4}[\Psi_{x}^{2}(1)
+\Psi_{x}^{2}(0)],$$
$$\frac{K\rho_3}{\rho_1}[\Phi_{x}\theta_{t}]_{x=0}^{x=1}\leq\frac{\varepsilon_4\rho_3K^{2}}{\rho_1^{2}b}[\Phi_{x}^{2}(1)+\Phi_{x}^{2}(0)]
+\frac{b\rho_3}{4\varepsilon_4}[\theta_{t}^{2}(1)+\theta_{t}^{2}(0)]$$
and
\bes
&&-\int_{0}^{1}\left(g*\theta_{x}\right)_{t}\Phi_{x}{\rm d}x=-\int_{0}^{1}\left(g(0)\theta_{x}+g'*\theta_{x}\right)\Phi_{x}{\rm d}x=-\int_{0}^{1}\left(g(t)\theta_{x}-g'\diamond\theta_{x}\right)\Phi_{x}{\rm d}x\nonumber\\
&\leq&-\frac{g(0)}{2\varepsilon_4}\int_{0}^{1}\left(g'\circ \theta_{x}\right){\rm d}x+\frac{\varepsilon_4}{2}\int_{0}^{1}\Phi_{x}^{2}{\rm d}x+\frac{\varepsilon_4}{2}\int_{0}^{1}\Phi_{x}^{2}{\rm d}x+\frac{g^{2}(0)}{2\varepsilon_4}\int_{0}^{1}\theta_{x}^{2}{\rm d}x.
\ees
Combining all the above estimates, we get the desired results.\quad \quad $\Box$

In order to absorb the boundary terms, appearing in \eqref{3.0}, we exploit as in \cite{riverra2002}, the following function:
$$q(x)=2-4x, \ \ x\in[0,1].$$
We will also introduce the functionals $J_{1} $ and $J_{2} $ defined by
\bes
J_{1}(t)=\rho_1\int_{0}^{1}\Phi_{t}q\Phi_{x}{\rm d}x
\ees
and
\bes
J_{2}(t)=\gamma\rho_2b\int_{0}^{1}\Psi_{t}q\Psi_{x}{\rm d}x+\frac{\beta b}{\delta}\int_{0}^{1}\left(\rho_3\theta_{t}+\gamma\Psi_{x}\right)q\left(\delta\theta_{x}-\left(g*\theta_{x}\right)\right){\rm d}x.
\ees
\begin{lemma}\label{8}
Let $(\Phi,\Psi,\theta,z)$ be a solution of \eqref{1.5}-\eqref{12.5}. For any $\varepsilon_4>0, $ we have
\bes
J'_{1}(t)\leq-K[\Phi_{x}^{2}(1)+\Phi_{x}^{2}(0)]+2\rho_1\int_{0}^{1}\Phi_{t}^{2}{\rm d}x+3K\int_{0}^{1}\Phi_{x}^{2}{\rm d}x+K\int_{0}^{1}\Psi_{x}^{2}{\rm d}x
\lbl{3.2}
\ees
and
\bes
J'_{2}(t)&\leq&-b^{2}\gamma[\Psi_{x}^{2}(1)+\Psi_{x}^{2}(0)]-\beta b\rho_3[\theta_{t}^{2}(1)+\theta_{t}^{2}(0)]+\varepsilon_{4}^{2}K^{2}\int_{0}^{1}|\Phi_{x}+\Psi|^{2}{\rm d}x+2\gamma\rho_2b\int_{0}^{1}\Psi_{t}^{2}{\rm d}x\nonumber\\
&&+\left(2b^{2}\gamma+\frac{\gamma^{2}b^{2}}{4\varepsilon_{4}^{2}}+\varepsilon_{4}\right)\int_{0}^{1}\Psi_{x}^{2}{\rm d}x+\left(2\beta b\rho_3+\frac{5\varepsilon_{4}}{2}\right)\int_{0}^{1}\theta_{t}^{2}{\rm d}x+C_{1}(\varepsilon_4)\int_{0}^{1}\theta_{x}^{2}{\rm d}x\nonumber\\
&&+\frac{3\varepsilon_{4}}{2}\int_{0}^{1}z^{2}(x,1,t){\rm d}x+\frac{2\beta b\int_{0}^{t}g(s){\rm d}s}{\delta^{2}\varepsilon_4}\left[4\delta\varepsilon_4+\beta b\mu_1^{2}+\beta b\mu_2^{2}\right]\int_{0}^{1}\left(g\circ\theta_{x}\right){\rm d}x\nonumber\\
&&-\frac{2g(0)\beta^{2}b^{2}\left(\gamma^{2}+\rho_{3}^{2}\right)}{\varepsilon_4\delta^{2}}\int_{0}^{1}\left(g'\circ\theta_{x}\right){\rm d}x,
\lbl{3.3}
\ees
where
$$C_{1}(\varepsilon_4)=\frac{4\beta b}{\delta}\left[2\left(\int_{0}^{t}g(s){\rm d}s\right)^{2}+\delta^{2}\right]+\frac{ 2\beta^{2}b^{2}}{\varepsilon_4}\left[\frac{g^{2}(0)(\gamma^2+\rho_3^{2})}{\delta^{2}}+\left(\mu_1^2+\mu_2^2\right)\left(1+\frac{\left(\int_{0}^{t}g(s){\rm d}s\right)^{2}}{\delta^{2}}\right)\right].
$$
\end{lemma}
{\bf Proof. } A direct differentiation of $J_{1} $ yields
\bes
J'_{1}(t)&=&K\int_{0}^{1}(\Phi_{x}+\Psi)_{x}q\Phi_{x}{\rm d}x+\rho_{1}\int_{0}^{1}\Phi_{t}q\Phi_{tx}{\rm d}x\nonumber\\
&=&K\int_{0}^{1}\Phi_{xx}q\Phi_{x}{\rm d}x+K\int_{0}^{1}\Psi_{x}q\Phi_{x}{\rm d}x-\frac{\rho_1}{2}\int_{0}^{1}q_{x}\Phi_{t}^{2}{\rm d}x\nonumber\\
&=&\frac{K}{2}[q\Phi_{x}^{2}]_{x=0}^{x=1}-\frac{K}{2}\int_{0}^{1}q_{x}\Phi_{x}^{2}{\rm d}x+K\int_{0}^{1}\Psi_{x}q\Phi_{x}{\rm d}x+2\rho_1\int_{0}^{1}\Phi_{t}^{2}{\rm d}x\nonumber\\
&=&-K[\Phi_{x}^{2}(1)+\Phi_{x}^{2}(0)]+2K\int_{0}^{1}\Phi_{x}^{2}{\rm d}x+2\rho_1\int_{0}^{1}\Phi_{t}^{2}{\rm d}x+K\int_{0}^{1}\Psi_{x}q\Phi_{x}{\rm d}x.
\ees
The Young's inequality applied to the last term gives the result.

Differentiating $J_{2}(t) $ along solutions of \eqref{1.5}, we find
\bes
J'_{2}(t)&=&\gamma b\int_{0}^{1}\left[b\Psi_{xx}-K(\Phi_{x}+\Psi)-\beta\theta_{tx}\right]q\Psi_{x}{\rm d}x+\gamma\rho_2b\int_{0}^{1}\Psi_{t}q\Psi_{tx}{\rm d}x\nonumber\\
&&+\frac{\beta b}{\delta}\int_{0}^{1}\left[(\delta\theta_{xx}-\gamma\Psi_{tx}-\mu_1\theta_{t}-\mu_2z(x,1,t)-g*\theta_{xx})+\gamma\Psi_{tx}\right]q(\delta\theta_{x}-g*\theta_{x}){\rm d}x\nonumber\\
&&+\frac{\beta b}{\delta}\int_{0}^{1}(\rho_{3}\theta_{t}+\gamma\Psi_{x})q\left(\delta\theta_{tx}-(g*\theta_{x})_{t}\right){\rm d}x.
\ees
By integration by part, we obtain
\bes
J'_{2}(t)&=&\frac{b^{2}\gamma}{2}[q\Psi_{x}^{2}]_{x=0}^{x=1}-\frac{b^{2}\gamma}{2}\int_{0}^{1}q_{x}\Psi_{x}^{2}{\rm d}x-K\gamma b\int_{0}^{1}(\Phi_{x}+\Psi)q\Psi_{x}{\rm d}x\nonumber\\
&&-\beta\gamma b\int_{0}^{1}\theta_{tx}q\Psi_{x}{\rm d}x-\frac{\gamma\rho_2 b}{2}\int_{0}^{1}q_{x}\Psi_{t}^{2}{\rm d}x-\frac{\beta b}{2\delta}\int_{0}^{1}q_{x}\left(\delta\theta_{x}-(g*\theta_{x})\right)^{2}{\rm d}x\nonumber\\
&&-\frac{\beta b \mu_1}{\delta}\int_{0}^{1}\theta_{t}q\left(\delta\theta_{x}-(g*\theta_{x})\right){\rm d}x-\frac{\beta b \mu_2}{\delta}\int_{0}^{1}z(x,1,t)q\left(\delta\theta_{x}-(g*\theta_{x})\right){\rm d}x+\frac{\beta b\rho_{3}}{2}[q\theta_{t}^{2}]_{x=0}^{x=1}\nonumber\\
&&-\frac{\beta b\rho_{3}}{2}\int_{0}^{1}q_{x}\theta_{t}^{2}{\rm d}x-\frac{\beta b\rho_{3}}{\delta}\int_{0}^{1}\theta_{t}q(g*\theta_{x})_{t}{\rm d}x+\gamma\beta b\int_{0}^{1}\Psi_{x}q\theta_{tx}{\rm d}x-\frac{b \beta \gamma }{\delta}\int_{0}^{1}\Psi_{x}q(g*\theta_{x})_{t}{\rm d}x\nonumber\\
&=&-b^{2}\gamma[\Psi_{x}^{2}(1)+\Psi_{x}^{2}(0)]+2b^{2}\gamma\int_{0}^{1}\Psi_{x}^{2}{\rm d}x-K\gamma b\int_{0}^{1}(\Phi_{x}+\Psi)q\Psi_{x}{\rm d}x+2\gamma\rho_2 b\int_{0}^{1}\Psi_{t}^{2}{\rm d}x\nonumber\\
&&-\frac{\beta b\gamma}{\delta}\int_{0}^{1}\Psi_{x}q(g*\theta_{x})_{t}{\rm d}x-\beta b\mu_1\int_{0}^{1}\theta_{t}q\theta_{x}{\rm d}x+\frac{\beta b\mu_1}{\delta}\int_{0}^{1}\theta_{t}q(g*\theta_{x}){\rm d}x+2\beta b\rho_{3}\int_{0}^{1}\theta_{t}^{2}{\rm d}x\nonumber\\
&&-\beta b\mu_2\int_{0}^{1}z(x,1,t)q\theta_{x}{\rm d}x+\frac{\beta b\mu_2}{\delta}\int_{0}^{1}z(x,1,t)q(g*\theta_{x}){\rm d}x-\beta b \rho_{3}[\theta_{t}^{2}(1)+\theta_{t}^{2}(0)]\nonumber\\
&&-\frac{\beta b\rho_{3}}{\delta}\int_{0}^{1}\theta_{t}q(g*\theta_{x})_{t}{\rm d}x+\frac{2\beta b}{\delta}\int_{0}^{1}\left(\delta\theta_{x}-(g*\theta_{x})\right)^{2}{\rm d}x.
\ees
Note that
\bes
&&-\frac{\beta b\gamma}{\delta}\int_{0}^{1}\Psi_{x}q(g*\theta_{x})_{t}{\rm d}x=-\frac{\beta b\gamma}{\delta}\int_{0}^{1}\Psi_{x}q(g(t)\theta_{x}-g'\diamond\theta_{x}){\rm d}x\nonumber\\
&\leq&\frac{2}{\varepsilon_{4}}\left(\frac{b\beta\gamma}{\delta}\right)^{2}g^{2}(0)\int_{0}^{1}\theta_{x}^{2}{\rm d}x+\varepsilon_{4}\int_{0}^{1}\Psi_{x}^{2}{\rm d}x-\frac{2}{\varepsilon_{4}}\left(\frac{b\beta\gamma}{\delta}\right)^{2}g(0)\int_{0}^{1}(g'\circ\theta_{x}){\rm d}x
\ees
and similarly
$$-\frac{\beta b\rho_3}{\delta}\int_{0}^{1}\theta_{t}q(g*\theta_{x})_{t}{\rm d}x\leq\frac{2}{\varepsilon_{4}}\left(\frac{b\beta\rho_3}{\delta}\right)^{2}g^{2}(0)\int_{0}^{1}\theta_{x}^{2}{\rm d}x+\varepsilon_{4}\int_{0}^{1}\theta_{t}^{2}{\rm d}x-\frac{2}{\varepsilon_{4}}\left(\frac{b\beta\rho_3}{\delta}\right)^{2}g(0)\int_{0}^{1}(g'\circ\theta_{x}){\rm d}x.$$
Moreover,
$$-\beta b\mu_1\int_{0}^{1}\theta_{t}q\theta_{x}{\rm d}x\leq\frac{\varepsilon_4}{2}\int_{0}^{1}\theta_{t}^{2}{\rm d}x+\frac{2\beta^{2}b^{2}\mu_1^{2}}{\varepsilon_4}\int_{0}^{1}\theta_{x}^{2}{\rm d}x,$$
$$-\beta b\mu_2\int_{0}^{1}z(x,1,t)q\theta_{x}{\rm d}x\leq\frac{\varepsilon_4}{2}\int_{0}^{1}z^{2}(x,1,t){\rm d}x+\frac{2\beta^{2}b^{2}\mu_2^{2}}{\varepsilon_4}\int_{0}^{1}\theta_{x}^{2}{\rm d}x,$$
$$-K\gamma b\int_{0}^{1}(\Phi_{x}+\Psi)q\Psi_{x}{\rm d}x\leq\varepsilon_{4}^{2}K^{2}\int_{0}^{1}|\Phi_{x}+\Psi|^{2}{\rm d}x+\frac{\gamma^{2}b^{2}}{4\varepsilon_{4}^{2}}\int_{0}^{1}\Psi_{x}^{2}{\rm d}x,$$
\bes
\frac{2\beta b}{\delta}\int_{0}^{1}\left(\delta\theta_{x}-(g*\theta_{x})\right)^{2}{\rm d}x&\leq&\frac{8\beta b}{\delta}\left(\int_{0}^{t}g(s){\rm d}s\right)^{2}\int_{0}^{1}\theta_{x}^{2}{\rm d}x+4\beta b\delta\int_{0}^{1}\theta_{x}^{2}{\rm d}x\nonumber\\
&&+\frac{8\beta b}{\delta}\left(\int_{0}^{t}g(s){\rm d}s\right)\int_{0}^{1}(g\circ\theta_{x}){\rm d}x
\ees
and
\bes
\frac{\beta b\mu_1}{\delta}\int_{0}^{1}\theta_{t}q(g*\theta_{x}){\rm d}x&\leq&\frac{\varepsilon_4}{2}\int_{0}^{1}\theta_{t}^{2}{\rm d}x+\frac{2}{\varepsilon_4}\left(\frac{\beta b\mu_1}{\delta}\right)^{2}\int_{0}^{t}g(s){\rm d}s\int_{0}^{1}(g\circ\theta_{x}){\rm d}x\nonumber\\
&&+\frac{\varepsilon_4}{2}\int_{0}^{1}\theta_{t}^{2}{\rm d}x+\frac{2}{\varepsilon_4}\left(\frac{\beta b\mu_1}{\delta}\right)^{2}\left(\int_{0}^{t}g(s){\rm d}s\right)^{2}\int_{0}^{1}\theta_{x}^{2}{\rm d}x,
\ees
\bes
\frac{\beta b\mu_2}{\delta}\int_{0}^{1}z(x,1,t)q(g*\theta_{x}){\rm d}x&\leq&\varepsilon_4\int_{0}^{1}z^{2}(x,1,t){\rm d}x+\frac{2}{\varepsilon_4}\left(\frac{\beta b\mu_2}{\delta}\right)^{2}\int_{0}^{t}g(s){\rm d}s\int_{0}^{1}(g\circ\theta_{x}){\rm d}x\nonumber\\
&&+\frac{2}{\varepsilon_4}\left(\frac{\beta b\mu_2}{\delta}\right)^{2}\left(\int_{0}^{t}g(s){\rm d}s\right)^{2}\int_{0}^{1}\theta_{x}^{2}{\rm d}x.
\ees
This completes the proof of the lemma.\quad \quad $\Box$

Now, let us define the following functional
$$I_{5}(t)=\rho_2\rho_3\int_{0}^{1}\int_{0}^{x}\theta_{t}(t,y){\rm d}y\Psi_{t}{\rm d}x.$$
\begin{lemma}\label{9}
Let $(\Phi,\Psi,\theta,z)$ be a solution of \eqref{1.5}-\eqref{12.5}. For any $\varepsilon_2>0 $ and $\eta_1>0, $ we have
\bes
I'_{5}(t)&\leq&-\frac{\rho_2\gamma}{4}\int_{0}^{1}\Psi_{t}^{2}{\rm d}x+\frac{\rho_{2}}{\gamma}\left[\delta^{2}+2\left(\int_{0}^{t}g(s){\rm d}s\right)^{2}\right]\int_{0}^{1}\theta_{x}^{2}{\rm d}x
+\frac{2\rho_2\mu_2^{2}}{\gamma}\int_{0}^{1}z^2(x,1,t){\rm d}x\nonumber\\
&&+\varepsilon_2(C_{p}+1)\int_{0}^{1}\Psi_{x}^{2}{\rm d}x+\frac{2\rho_2}{\gamma}\int_{0}^{t}g(s){\rm d}s\int_{0}^{1}(g\circ\theta_{x}){\rm d}x+\frac{b^2}{4\eta_1}\Psi_{x}^{2}(1)+\frac{\rho_3b}{4\eta_1}\theta_{t}^{2}(1)\nonumber\\
&&+\left[\beta\rho_3+\frac{\rho_2\mu_1^{2}}{\gamma}+\frac{\rho_3^{2}}{4\varepsilon_2}\left(2K^{2}+b^{2}\right)+\eta_1\left(\rho_{3}^{2}+\frac{\rho_3\beta^{2}}{b}\right)\right]\int_{0}^{1}\theta_{t}^{2}{\rm d}x\nonumber\\
&&+\varepsilon_2C_{p}\int_{0}^{1}\Phi_{x}^{2}{\rm d}x.
\lbl{3.5}
\ees
\end{lemma}
{\bf Proof. } Differentiating the functional $I_{5} $ and using Eqs. $\eqref{1.5},$ we obtain
\bes
I'_{5}(t)&=&\rho_2\int_{0}^{1}\int_{0}^{x}[\delta\theta_{xx}-\gamma\Psi_{tx}-\mu_1\theta_{t}-\mu_2z(y,1,t)-g*\theta_{xx}]{\rm d}y\Psi_{t}{\rm d}x\nonumber\\
&&+\rho_{3}\int_{0}^{1}\int_{0}^{x}\theta_{t}{\rm d}y[b\Psi_{xx}-K(\Phi_{x}+\Psi)-\beta\theta_{tx}]{\rm d}x\nonumber\\
&=&-\gamma\rho_2\int_{0}^{1}\Psi_{t}^{2}{\rm d}x+\rho_{2}\delta\int_{0}^{1}\theta_{x}\Psi_{t}{\rm d}x-\rho_{2}\int_{0}^{1}(g*\theta_{x})\Psi_{t}{\rm d}x-\mu_1\rho_{2}\int_{0}^{1}\int_{0}^{x}\theta_{t}{\rm d}y\Psi_{t}{\rm d}x\nonumber\\
&&-\mu_2\rho_{2}\int_{0}^{1}\int_{0}^{x}z(y,1,t){\rm d}y\Psi_{t}{\rm d}x-K\rho_{3}\int_{0}^{1}\int_{0}^{x}\theta_{t}{\rm d}y\Psi{\rm d}x-\rho_{3}b\int_{0}^{1}\theta_{t}\Psi_{x}{\rm d}x\nonumber\\
&&+K\rho_{3}\int_{0}^{1}\theta_{t}\Phi{\rm d}x+\beta\rho_{3}\int_{0}^{1}\theta_{t}^{2}{\rm d}x+\rho_{3}\left[\int_{0}^{x}\theta_{t}{\rm d}y(b\Psi_{x}-\beta\theta_{t})\right]_{x=0}^{x=1}.
\ees
Using Young's and Poincar$\acute{e}$'s inequalities, we find
\bes
I'_{5}(t)&\leq&-\frac{\gamma\rho_2}{2}\int_{0}^{1}\Psi_{t}^{2}{\rm d}x+\frac{\rho_{2}\delta^{2}}{\gamma}\int_{0}^{1}\theta_{x}^{2}{\rm d}x+\varepsilon_{2}(1+C_{p})\int_{0}^{1}\Psi_{x}^{2}{\rm d}x+\varepsilon_{2}C_{p}\int_{0}^{1}\Phi_{x}^{2}{\rm d}x\nonumber\\
&&+\left[\beta\rho_{3}+\frac{2\rho_2\mu_1^{2}}{\gamma}+\frac{\rho_3^{2}}{4\varepsilon_2}\left(2K^{2}+b^{2}\right)\right]\int_{0}^{1}\theta_{t}^{2}{\rm d}x+\frac{2\rho_2\mu_2^2}{\gamma}\int_{0}^{1}z^2(x,1,t){\rm d}x\nonumber\\
&&-\rho_{2}\int_{0}^{1}(g*\theta_{x})\Psi_{t}{\rm d}x+\rho_{3}\left[\int_{0}^{x}\theta_{t}{\rm d}y(b\Psi_{x}-\beta\theta_{t})\right]_{x=0}^{x=1}.
\ees
We can estimate terms in the right side as follows
\bes
&&\rho_{2}\int_{0}^{1}(g*\theta_{x})\Psi_{t}{\rm d}x=-\rho_{2}\int_{0}^{1}(g\diamond\theta_{x})\Psi_{t}{\rm d}x+\rho_{2}\int_{0}^{t}g(s){\rm d}s\int_{0}^{1}\theta_{x}\Psi_{t}{\rm d}x\nonumber\\
&\leq&\frac{\rho_{2}\gamma}{4}\int_{0}^{1}\Psi_{t}^{2}{\rm d}x+\frac{2\rho_{2}}{\gamma}\int_{0}^{t}g(s){\rm d}s\int_{0}^{1}(g\circ\theta_{x}){\rm d}x+\frac{2\rho_{2}}{\gamma}\left(\int_{0}^{t}g(s){\rm d}s\right)^{2}\int_{0}^{1}\theta_{x}^{2}{\rm d}x
\ees
and
$$\rho_{3}\left[\int_{0}^{x}\theta_{t}^{2}{\rm d}y(b\Psi_{x}-\beta\theta_{t})\right]_{x=0}^{x=1}\leq\frac{b^{2}}{4\eta_{1}}\Psi_{x}^{2}(1)+\frac{\rho_{3}b}{4\eta_{1}}\theta_{t}^{2}(1)
+\eta_{1}\left(\rho_{3}^{2}+\frac{\rho_{3}\beta^{2}}{b}\right)\int_{0}^{1}\theta_{t}^{2}{\rm d}x.$$
The proof is completed.\quad \quad $\Box$

Finally, as in \cite{kirane2011}, we introduce the functional
$$I_{6}(t)=\int_{0}^{1}\int_{0}^{1}e^{-2\tau\rho}z^{2}(x,\rho,t){\rm d}\rho{\rm d}x.$$
Then the following result holds:
\begin{lemma}\label{10}
Let $(\Phi,\Psi,\theta,z)$ be a solution of \eqref{1.5}-\eqref{12.5}. Then, we have
\bes
I'_{6}(t)\leq-2I_{6}(t)-\frac{c}{\tau}\int_{0}^{1}z^2(x,1,t){\rm d}x+\frac{1}{\tau}\int_{0}^{1}\theta^2_{t}{\rm d}x.
\lbl{3.6}
\ees
where $c $ is a positive constant.
\end{lemma}
{\bf Proof. } Differentiating the functional $I_{6}, $ we have
\bes
I'_{6}(t)&=&-\frac{2}{\tau}\int_{0}^{1}\int_{0}^{1}e^{-2\tau\rho}z z_{\rho}(x,\rho,t){\rm d}\rho{\rm d}x\nonumber\\
&=&-2\int_{0}^{1}\int_{0}^{1}e^{-2\tau\rho}z^2(x,\rho,t){\rm d}\rho{\rm d}x-\frac{1}{\tau}\int_{0}^{1}\int_{0}^{1}\frac{\partial}{\partial\rho}\left(e^{-2\tau\rho}z^2(x,\rho,t)\right){\rm d}\rho{\rm d}x\nonumber\\
&=&-2I_6(t)+\frac{1}{\tau}\int_{0}^{1}\theta^2_{t}{\rm d}x-\frac{1}{\tau}\int_{0}^{1}e^{-2\tau}z^2(x,1,t){\rm d}x.
\ees
The above equality implies that there exists a positive constant $c $ such that \eqref{3.6} holds.\quad \quad $\Box$

Now, we define the Lyapunov functional $\mathcal{L}$ as follows
\bes
\mathcal{L}(t)&=&NE(t)+N_1I_1(t)+N_2I_2(t)+\frac{\upsilon}{4}I_3(t)+\upsilon I_4(t)+N_5I_5(t)+I_6(t)\nonumber\\
&&+\upsilon\varepsilon_{4}\left(\frac{1}{K}+\frac{\rho_{3}K}{\rho_{1}^{2}b}\right)J_{1}(t)+\frac{1}{2\varepsilon_{4}}J_{2}(t),\ \ t\geq0,
\lbl{3.7}
\ees
where $N, N_1, N_2, N_5 $ are positive constants to be chosen properly later and $\upsilon=\min\{\gamma, \beta\}.$ For large $N, $ we can verify that, for some $m, M>0,$
\bes
mE(t)\leq\mathcal{L}(t)\leq ME(t), \ \ t\geq0.
\lbl{3.8}
\ees

Taking into account \eqref{1.8}, \eqref{2.7}, \eqref{2.8}, \eqref{2.9}, \eqref{3.0}, \eqref{3.2}, \eqref{3.3}, \eqref{3.5}, \eqref{3.6} and the relations
\bes
\int_{0}^{1}\Phi_{x}^{2}{\rm d}x\leq2\int_{0}^{1}|\Phi_{x}+\Psi|^{2}{\rm d}x+2C_{p}\int_{0}^{1}\Psi_{x}^{2}{\rm d}x,
\lbl{5.5}
\ees
we arrive at
\bes
\mathcal{L'}(t)&\leq&-\left\{\frac{\rho_1\upsilon}{4}-N_{1}\varepsilon_{1}-2\varepsilon_{4}\upsilon\left[1+\rho_{1}\left(\frac{1}{K}+\frac{\rho_{3}K}{\rho_{1}^{2}b}\right)\right] \right\}\int_{0}^{1}\Phi_{t}^{2}{\rm d}x\nonumber\\
&&-\left\{\frac{\rho_{2}\gamma N_{5}}{4}-N_{1}\left(\frac{3\rho_{2}}{2}+\frac{\rho_{1}^{2}C_{p}}{4\varepsilon_{1}}\right)-\frac{3\rho_{2}\upsilon}{4}-\frac{\gamma\rho_{2}b}{\varepsilon_{4}} \right\}\int_{0}^{1}\Psi_{t}^{2}{\rm d}x\nonumber\\
&&-\left\{\frac{3\upsilon K}{4}-2\varepsilon_2C_{p}N_5-\varepsilon_{4}\left[6\upsilon K\left(\frac{1}{K}+\frac{\rho_{3}K}{\rho_{1}^{2}b}\right)+2\upsilon+\frac{K^{2}}{2}\right]\right\}\int_{0}^{1}|\Phi_{x}+\Psi|^{2}{\rm d}x\nonumber\\
&&-\left\{N_{1}b-\varepsilon_{2}\left(N_{2}+N_{5}\left(1+C_{p}+2C_{p}^{2}\right)\right)-\frac{1}{2\varepsilon_{4}}\left[2b^{2}\gamma+\frac{\gamma^{2}b^{2}}{4\varepsilon_{4}^{2}}+\varepsilon_{4}\right]\right.\nonumber\\
&&\left.-\varepsilon_{4}\left(K\upsilon+6C_{p}K\upsilon\right)\left(\frac{1}{K}+\frac{\rho_{3}K}{\rho_{1}^{2}b}\right)-\frac{3b\upsilon}{8}-\varepsilon_{4}\left(\upsilon+2\upsilon C_{p}\right)\right\}\int_{0}^{1}\Psi_{x}^{2}{\rm d}x\nonumber\\
&&-\left\{Nm_0-N_{2}\left(\rho_{3}+\frac{\gamma^{2}}{4\varepsilon_{2}}\right)-\left[\frac{\upsilon}{4\varepsilon_{4}}(\delta^{2}+\mu_1^{2})+\frac{\beta b\rho_{3}}{\varepsilon_{4}}+\frac{5}{4}\right]-\frac{\beta^2\upsilon}{8b}\right.\nonumber\\
&&\left.-N_{5}\left[\beta\rho_{3}+\frac{\rho_2\mu_1^{2}}{\gamma}+\frac{\rho_3^{2}}{4\varepsilon_{2}}\left(2K^{2}+b^{2}\right)+\eta_{1}\left(\rho_{3}^{2}+\frac{\rho_{3}\beta^{2}}{b}\right)\right]-\frac{1}{\tau}
\right\}\int_{0}^{1}\theta_{t}^{2}{\rm d}x\nonumber\\
&&-\left\{\frac{\lambda N_{2}}{2}-\frac{N_{1}\beta^{2}}{2\rho_{2}}-\frac{1}{2\varepsilon_{4}}\left(\upsilon g^{2}(0)+C_{1}(\varepsilon_{4})\right)-\frac{N_{5}\rho_{2}}{\gamma}\left(\delta^{2}+2\overline{g}^{2}\right)\right\}\int_{0}^{1}\theta_{x}^{2}{\rm d}x\nonumber\\
&&+\left\{\frac{N_{2}\overline{g}}{\lambda}+\frac{\beta b\overline{g}}{\varepsilon_4^2\delta^2}\left(4\delta\varepsilon_4+\beta b\mu_1^2+\beta b \mu_2^2\right)+\frac{2\rho_2\overline{g}N_5}{\gamma}\right\}\int_{0}^{1}\left(g\circ\theta_{x}\right){\rm d}x\nonumber\\
&&+\left\{\frac{\beta N}{2}-\frac{\upsilon g(0)}{2\varepsilon_{4}}-\frac{g(0)b^{2}\beta^{2}\left(\gamma^{2}+\rho_{3}^{2}\right)}{\varepsilon_{4}^{2}\delta^{2}} \right\}\int_{0}^{1}\left(g'\circ\theta_{x}\right){\rm d}x-2I_6(t)\nonumber\\
&&-\left\{Nm_0+\frac{c}{\tau}-\frac{\mu_2^2N_2C_p}{\lambda}-\frac{\mu_2^2\upsilon}{4\varepsilon_4}-\frac{3}{4}-\frac{2\rho_2\mu_2^2N_5}{\gamma}\right\}\int_{0}^{1}z^2(x,1,t){\rm d}x,
\lbl{3.9}
\ees
where $m_0=\min\left\{\beta\left(\mu_1-\frac{\xi}{2\tau}-\frac{\mu_2}{2}\right), \beta\left(\frac{\xi}{2\tau}-\frac{\mu_2}{2}\right)\right\}. $ At this point, we need to choose our constants very carefully. First, let us pick $\eta_{1}=\frac{N_5\varepsilon_4}{\upsilon}$ and choose
\bes
\varepsilon_{4}\leq\min\left\{\frac{\rho_{1}}{16}\left[1+\rho_{1}\left(\frac{1}{K}+\frac{\rho_{3}K}{\rho_{1}^{2}b}\right)\right]^{-1}, \frac{3\upsilon K}{8}\left[6\upsilon K\left(\frac{1}{K}+\frac{\rho_{3}K}{\rho_{1}^{2}b}\right)+2\upsilon+\frac{K^{2}}{2}\right]^{-1}\right\}.
\lbl{4.0}
\ees
Second, we select $N_{1} $ sufficiently large such that
\bes
\frac{N_{1}b}{2}&\geq&\frac{1}{2\varepsilon_{4}}\left[2b^{2}\gamma+\frac{\gamma^{2}b^{2}}{4\varepsilon_{4}^{2}}+\varepsilon_{4}\right]+\varepsilon_{4}K\upsilon\left(1+6C_{p}\right)\left(\frac{1}{K}+\frac{\rho_{3}K}{\rho_{1}^{2}b}\right)\nonumber\\&&+\frac{3b\upsilon}{8}+\varepsilon_{4}\upsilon\left(1+2 C_{p}\right),
\lbl{4.1}
\ees
then we choose $\varepsilon_{1} $ so small that
$$\varepsilon_{1}\leq\frac{\rho_{1}\upsilon}{16N_{1}}.$$
Next, we choose $N_{5} $ sufficiently large so that
$$N_{5}\geq\frac{8}{\rho_{2}\gamma}\left[N_{1}\left(\frac{3\rho_{2}}{2}+\frac{\rho_{1}^{2}C_{p}}{4\varepsilon_{1}}\right)+\frac{\gamma\rho_{2}b}{\varepsilon_{4}} +\frac{3\rho_{2}\upsilon}{4}\right],$$
and also select $N_{2} $ sufficiently large so that
$$\frac{\lambda N_{2}}{4}>\frac{N_{1}\beta^{2}}{2\rho_{2}}+\frac{1}{2\varepsilon_{4}}\left(\upsilon g^{2}(0)+C_{1}(\varepsilon_{4})\right)+\frac{N_{5}\rho_{2}}{\gamma}\left(\delta^{2}+2\overline{g}^{2}\right).$$
Furthermore, we select $\varepsilon_{2} $  satisfies
$$\varepsilon_{2}<\min\left\{\frac{N_{1}b}{2\left(N_{2}+N_{5}\left(1+C_{p}+2C_{p}^{2}\right)\right)}, \frac{3K\upsilon}{16N_{5}C_{p}}\right\}.$$
Finally, we choose $N $ large enough so that \eqref{3.8} remains valid and \eqref{3.9} takes the form
\bes
\mathcal{L'}(t)&\leq&-C_{1}\int_{0}^{1}\left(\Phi_{t}^{2}+\Psi_{t}^{2}+|\Phi_{x}+\Psi|^{2}+\Psi_{x}^{2}+\theta_{t}^{2}+\theta_{x}^{2}+\int_{0}^{1}z^2(x,\rho,t){\rm d}\rho\right){\rm d}x+C_{2}\int_{0}^{1}(g\circ\theta_{x}){\rm d}x\nonumber\\
&\leq&-CE(t)+C_{3}\int_{0}^{1}(g\circ\theta_{x}){\rm d}x,
\lbl{4.2}
\ees
where $C_{1}, C_{2}, C_{3}, $ and $C $ are positive constants.

{\large\textbf{4.2 The case $\mu_2=\mu_1$}}

If $\mu_1=\mu_2=\mu,$ then we can choose $\xi=\tau\mu$ in \eqref{1.6} and Lemma \ref{2} takes the form
\begin{lemma}\lbl{11}
Let $(\Phi,\Psi,\theta,z)$ be a solution of \eqref{1.5}-\eqref{12.5}. Assume that $\mu_1=\mu_2=\mu $ and $g $ satisfies $(H1) $ and $(H2), $ $\xi=\tau\mu. $ Then, the energy functional defined by \eqref{1.7} is a non-increasing function and it satisfies
\bes
E'(t)&\leq&-\frac{\beta}{2}g(t)\int_{0}^{1}\theta_{x}^{2}{\rm d}x+\frac{\beta}{2}\int_{0}^{1}(g'\circ\theta_{x}){\rm d}x\leq 0, \ \ t\geq0.
\lbl{4.7}
\ees
\end{lemma}
The proof of Lemma \ref{11} is an immediate consequence of Lemma \ref{2}, by choosing $\xi=\tau\mu.$

If $\mu_1=\mu_2=\mu, $ we need some additional negative term of $\int_{0}^{1}\theta_{t}^{2}{\rm d}x. $ For this purpose, let us introduce the functional
$$I_{7}(t)=-\rho_3\int_{0}^{1}\theta_{t}\left(g\diamond\theta\right){\rm d}x.$$
Then, we have the following estimate:
\begin{lemma}\label{13}
Let $(\Phi,\Psi,\theta,z)$ be a solution of \eqref{1.5}-\eqref{12.5}. Then for any $\varepsilon_{7}>0 $ and $\eta_2>0, $ we have
\bes
I'_{7}(t)&\leq&-\left(\rho_3\int_{0}^{t}g(s){\rm d}s-\eta_2\right)\int_{0}^{1}\theta_{t}^{2}{\rm d}x+\varepsilon_{7}\int_{0}^{1}\Psi_{t}^{2}{\rm d}x+\varepsilon_{7}\left[1+\left(\int_{0}^{t}g(s){\rm d}s\right)^{2}\right]\int_{0}^{1}\theta_{x}^{2}{\rm d}x\nonumber\\
&&+\varepsilon_7\int_{0}^{1}z^2(x,1,t){\rm d}x+C_2(\varepsilon_7)\int_{0}^{1}\left(g\circ\theta_{x}\right){\rm d}x-\frac{\rho_3^{2}}{2\eta_2}g(0)C_{p}\int_{0}^{1}\left(g'\circ\theta_{x}\right){\rm d}x,
\lbl{4.9}
\ees
where $$C_2(\varepsilon_7)=\frac{\int_{0}^{t}g(s){\rm d}s}{4\varepsilon_7}\left(\delta^{2}+\gamma^{2}+4\varepsilon_{7}^{2}+2+\mu_2^2C_{p}\right)
+\frac{\mu_1^{2}C_{p}}{2\eta_2}\int_{0}^{t}g(s){\rm d}s.$$
\end{lemma}
{\bf Proof. } A simple differentiation leads to
\bes
I'_{7}(t)&=&-\rho_3\int_{0}^{1}\theta_{t}\left(g\diamond\theta\right)_{t}{\rm d}x-\rho_3\int_{0}^{1}\theta_{tt}\left(g\diamond\theta\right){\rm d}x\nonumber\\
&=&-\left(\rho_3\int_{0}^{t}g(s){\rm d}s\right)\int_{0}^{1}\theta_{t}^{2}{\rm d}x-\rho_3\int_{0}^{1}\theta_{t}\left(g'\diamond\theta\right){\rm d}x-\int_{0}^{1}\left(\int_{0}^{t}g(t-s)\theta_{x}(s){\rm d}s\right)\left(g\diamond\theta_{x}\right){\rm d}x\nonumber\\
&&+\delta\int_{0}^{1}\theta_{x}\left(g\diamond\theta_{x}\right){\rm d}x-\gamma\int_{0}^{1}\Psi_{t}\left(g\diamond\theta_{x}\right){\rm d}x+\mu_1\int_{0}^{1}\theta_{t}\left(g\diamond\theta\right){\rm d}x\nonumber\\
&&+\mu_2\int_{0}^{1}z(x,1,t)\left(g\diamond\theta\right){\rm d}x.
\lbl{5.0}
\ees
Terms in the right side of \eqref{5.0} are estimated as follows. Using Young's inequalities and Lemma \ref{1}, we obtain, for all $\eta_{2}>0,$
\bes
-\rho_3\int_{0}^{1}\theta_{t}\left(g'\diamond\theta\right){\rm d}x&\leq&\frac{\eta_2}{2}\int_{0}^{1}\theta_{t}^{2}{\rm d}x+\frac{\rho_{3}^{2}}{2\eta_2}\int_{0}^{t}(-g'(s){\rm d}s)\int_{0}^{1}(-g'\circ\theta){\rm d}x\nonumber\\
&\leq&\frac{\eta_2}{2}\int_{0}^{1}\theta_{t}^{2}{\rm d}x-\frac{\rho_{3}^{2}}{2\eta_2}g(0)C_{p}\int_{0}^{1}(g'\circ\theta_{x}){\rm d}x.
\ees
Similarly, for any $\varepsilon_7>0,$ we have
$$\delta\int_{0}^{1}\theta_{x}\left(g\diamond\theta_{x}\right){\rm d}x\leq\varepsilon_7\int_{0}^{1}\theta_{x}^{2}{\rm d}x+\frac{\delta^{2}}{4\varepsilon_7}\int_{0}^{t}g(s){\rm d}s\int_{0}^{1}\left(g\circ\theta_{x}\right){\rm d}x,$$
$$-\gamma\int_{0}^{1}\Psi_{t}\left(g\diamond\theta_{x}\right){\rm d}x\leq\varepsilon_7\int_{0}^{1}\Psi_{t}^{2}{\rm d}x+\frac{\gamma^{2}}{4\varepsilon_7}\int_{0}^{t}g(s){\rm d}s\int_{0}^{1}\left(g\circ\theta_{x}\right){\rm d}x,$$
$$\mu_1\int_{0}^{1}\theta_{t}\left(g\diamond\theta\right){\rm d}x\leq\frac{\eta_2}{2}\int_{0}^{1}\theta_{t}^{2}{\rm d}x+\frac{\mu_1^{2}C_{p}}{2\eta_2}\int_{0}^{t}g(s){\rm d}s\int_{0}^{1}\left(g\circ\theta_{x}\right){\rm d}x$$
and
$$\mu_2\int_{0}^{1}z(x,1,t)\left(g\diamond\theta\right){\rm d}x\leq\varepsilon_7\int_{0}^{1}z^{2}(x,1,t){\rm d}x+\frac{\mu_2^{2}C_{p}}{4\varepsilon_7}\int_{0}^{t}g(s){\rm d}s\int_{0}^{1}\left(g\circ\theta_{x}\right){\rm d}x.$$
Finally,
\bes
&&-\int_{0}^{1}\left(\int_{0}^{t}g(t-s)\theta_{x}(s){\rm d}s\right)\left(g\diamond\theta_{x}\right){\rm d}x\nonumber\\
&\leq&\frac{\varepsilon_7}{2}\int_{0}^{1}\left(\int_{0}^{t}g(t-s)\left(\theta_{x}(t)-\theta_{x}(x)-\theta_{x}(t)\right){\rm d}s\right)^{2}{\rm d}x+\frac{1}{2\varepsilon_7}\int_{0}^{1}\left(g\diamond\theta_{x}\right)^{2}{\rm d}x\nonumber\\
&\leq&\varepsilon_7\left(\int_{0}^{t}g(s){\rm d}s\right)^{2}\int_{0}^{1}\theta_{x}^{2}{\rm d}x+\left(\varepsilon_7+\frac{1}{2\varepsilon_7}\right)\int_{0}^{1}\left(g\diamond\theta_{x}\right)^{2}{\rm d}x\nonumber\\
&\leq&\varepsilon_7\left(\int_{0}^{t}g(s){\rm d}s\right)^{2}\int_{0}^{1}\theta_{x}^{2}{\rm d}x+\left(\varepsilon_7+\frac{1}{2\varepsilon_7}\right)\int_{0}^{t}g(s){\rm d}s\int_{0}^{1}\left(g\circ\theta_{x}\right){\rm d}x.
\ees
Therefore, the assertion of the lemma follows by combining all the above estimates.\quad \quad $\Box$

Now, we define the following Lyapunov functional $\mathcal{L} $ as:
\bes
\mathcal{L}(t)&=&NE(t)+N_1I_1(t)+N_2I_2(t)+\frac{\upsilon}{4}I_3(t)+\upsilon I_4(t)+N_5I_5(t)+N_6I_6(t)+N_7I_7(t)\nonumber\\
&&+\upsilon\varepsilon_{4}\left(\frac{1}{K}+\frac{\rho_{3}K}{\rho_{1}^{2}b}\right)J_{1}(t)+\frac{1}{2\varepsilon_{4}}J_{2}(t),\ \ t\geq0,
\lbl{5.6}
\ees
where $N, N_1, N_2, N_5, N_6, N_7 $ are positive real numbers which will be chosen later.

Since $g $ is continuous and $g(0)>0, $ then for any $t\geq t_{0}>0, $ we have
$$\int_{0}^{t}g(s){\rm d}s\geq\int_{0}^{t_{0}}g(s){\rm d}s=g_{0}.$$
Then, using the estimates \eqref{2.7}, \eqref{2.8}, \eqref{2.9}, \eqref{3.0}, \eqref{3.2}, \eqref{3.3}, \eqref{3.5}, \eqref{3.6}, \eqref{4.7}, \eqref{4.9} and algebraic inequality \eqref{5.5}, we get
\bes
\mathcal{L'}(t)&\leq&-\left\{\frac{\rho_1\upsilon}{4}-N_{1}\varepsilon_{1}-2\varepsilon_{4}\upsilon\left[1+\rho_{1}\left(\frac{1}{K}+\frac{\rho_{3}K}{\rho_{1}^{2}b}\right)\right] \right\}\int_{0}^{1}\Phi_{t}^{2}{\rm d}x-2N_6I_6(t)\nonumber\\
&&-\left\{\frac{\lambda N_{2}}{2}-\frac{N_{1}\beta^{2}}{2\rho_{2}}-N_7\varepsilon_7\left(1+\overline{g}^2\right)-\frac{1}{2\varepsilon_{4}}\left(\upsilon g^{2}(0)+C_{1}(\varepsilon_{4})\right)-\frac{N_{5}\rho_{2}}{\gamma}\left(\delta^{2}+2\overline{g}^{2}\right)\right\}\int_{0}^{1}\theta_{x}^{2}{\rm d}x\nonumber\\
&&-\left\{N_7\left(\rho_3g_0-\eta_2\right)-N_{2}\left(\rho_{3}+\frac{\gamma^{2}}{4\varepsilon_{2}}\right)-\left[\frac{\upsilon}{4\varepsilon_{4}}(\delta^{2}+\mu^{2})+\frac{\beta b\rho_{3}}{\varepsilon_{4}}+\frac{5}{4}\right]-\frac{\beta^2\upsilon}{8b}\right.\nonumber\\
&&\left.-N_{5}\left[\beta\rho_{3}+\frac{2\rho_2\mu^{2}}{\gamma}+\frac{\rho_3^{2}}{4\varepsilon_{2}}\left(2K^{2}+b^{2}\right)+\eta_{1}\left(\rho_{3}^{2}+\frac{\rho_{3}\beta^{2}}{b}\right)\right]-\frac{N_6}{\tau}
\right\}\int_{0}^{1}\theta_{t}^{2}{\rm d}x\nonumber\\
&&-\left\{\frac{3\upsilon K}{4}-2\varepsilon_2C_{p}N_5-\varepsilon_{4}\left[6\upsilon K\left(\frac{1}{K}+\frac{\rho_{3}K}{\rho_{1}^{2}b}\right)+2\upsilon+\frac{K^{2}}{2}\right]\right\}\int_{0}^{1}|\Phi_{x}+\Psi|^{2}{\rm d}x\nonumber\\
&&-\left\{\frac{\rho_{2}\gamma N_{5}}{4}-N_{1}\left(\frac{3\rho_{2}}{2}+\frac{\rho_{1}^{2}C_{p}}{4\varepsilon_{1}}\right)-N_7\varepsilon_7-\frac{3\rho_{2}\upsilon}{4}-\frac{\gamma\rho_{2}b}{\varepsilon_{4}} \right\}\int_{0}^{1}\Psi_{t}^{2}{\rm d}x\nonumber\\
&&-\left\{N_{1}b-\varepsilon_{2}\left(N_{2}+N_{5}\left(1+C_{p}+2C_{p}^{2}\right)\right)-\frac{1}{2\varepsilon_{4}}\left[2b^{2}\gamma+\frac{\gamma^{2}b^{2}}{4\varepsilon_{4}^{2}}+\varepsilon_{4}\right]\right.\nonumber\\
&&\left.-\varepsilon_{4}\left(K\upsilon+6C_{p}K\upsilon\right)\left(\frac{1}{K}+\frac{\rho_{3}K}{\rho_{1}^{2}b}\right)-\frac{3b\upsilon}{8}-\varepsilon_{4}\left(\upsilon+2\upsilon C_{p}\right)\right\}\int_{0}^{1}\Psi_{x}^{2}{\rm d}x\nonumber\\
&&+\left\{\frac{N_{2}\overline{g}}{2\lambda}+C_2(\varepsilon_7)+\frac{\beta b\overline{g}}{\varepsilon_4^2\delta^2}\left(4\delta\varepsilon_4+2\beta b\mu^2\right)+\frac{2\rho_2\overline{g}N_5}{\gamma}\right\}\int_{0}^{1}\left(g\circ\theta_{x}\right){\rm d}x\nonumber\\
&&+\left\{\frac{\beta N}{2}-\frac{\rho_3^2g(0)C_p}{2\eta_2}N_7-\frac{\upsilon g(0)}{2\varepsilon_{4}}-\frac{g(0)b^{2}\beta^{2}\left(\gamma^{2}+\rho_{3}^{2}\right)}{\varepsilon_{4}^{2}\delta^{2}} \right\}\int_{0}^{1}\left(g'\circ\theta_{x}\right){\rm d}x\nonumber\\
&&-\left\{\frac{N_6c}{\tau}-\frac{\mu^2N_2C_p}{\lambda}-N_7\varepsilon_7-\frac{\mu^2\upsilon}{4\varepsilon_4}-\frac{3}{4}-\frac{2\rho_2\mu^2N_5}{\gamma}\right\}\int_{0}^{1}z^2(x,1,t){\rm d}x.
\lbl{5.7}
\ees
Now, our goal is to choose our constants in \eqref{5.7} in order to get the negative coefficients on the right-hand side of \eqref{5.7}. To this end, let us pick $\eta_{1}=\frac{N_{5}\varepsilon_{4}}{\upsilon}, \eta_{2}=\frac{1}{4N_{7}} $ and we pick $\varepsilon_4, N_1, \varepsilon_1, N_5, N_2, \varepsilon_2 $ in the same order with the same values as the case $\mu_2<\mu_1$, respectively. Then we pick $N_6 $ large enough such that
$$\frac{N_6c}{2\tau}>\frac{N_2\mu^2C_p}{\lambda}+\frac{\mu^2\upsilon}{4\varepsilon_4}+\frac{3}{4}+\frac{2\rho_2\mu^2N_5}{\gamma}.$$
After that, we choose $N_7 $ sufficiently large so that
\bes
\frac{N_{7}\rho_{3}g_{0}}{2}-\frac{1}{8}&>&N_{2}\left(\rho_{3}+\frac{\gamma^{2}}{4\varepsilon_{2}}\right)+\frac{\upsilon\left(\delta^2+\mu^2\right)}{4\varepsilon_{4}}+\frac{\beta b\rho_{3}}{\varepsilon_{4}}+\frac{5}{4}+\frac{N_6}{\tau}+\frac{\beta^2\upsilon}{8b}\nonumber\\
&&+N_{5}\left[\beta\rho_{3}+\frac{2\rho_2\mu^{2}}{\gamma}+\frac{\rho_3^{2}}{4\varepsilon_{2}}\left(2K^{2}+b^{2}\right)+\eta_1\left(\rho_3^2+\frac{\rho_3\beta^2}{b}\right)\right].
\ees
Furthermore, choosing $\varepsilon_7 $ sufficiently small such that
$$\varepsilon_{7}<\min\left\{\frac{N_6c}{2\tau N_7}, \frac{\rho_{2}\gamma N_{5}}{8N_{7}}, \frac{\lambda N_{2}}{4N_{7}\left(1+\overline{g}^{2}\right)}\right\}.$$
Once all the above constants are fixed, we pick $N $ large enough such that there exists two positive constants $\hat{C} $ and $\hat{C_3} $
\bes
\mathcal{L'}(t)\leq-\hat{C}E(t)+\hat{C_3}\int_{0}^{1}(g\circ\theta_{x}){\rm d}x, \ \ t\geq t_0.
\lbl{80.888}
\ees

From \eqref{4.2} and \eqref{80.888}, we can know that the Lyapunov functionals $\mathcal{L}$ are of the same form under the two cases: $\mu_2<\mu_1 $ and $\mu_2=\mu_1 $
\bes
\mathcal{L'}(t)\leq-CE(t)+C_3\int_{0}^{1}(g\circ\theta_{x}){\rm d}x, \ \ t\geq t_0.
\lbl{88.88}
\ees

\noindent{\bf Continuity of the proof of Theorem \ref{3}.}  Multiplying \eqref{88.88} by $\zeta(t) $ gives
\bes
\zeta(t)\mathcal{L'}(t)\leq-C\zeta(t)E(t)+C_{3}\zeta(t)\int_{0}^{1}(g\circ\theta_{x}){\rm d}x.
\lbl{4.3}
\ees
The last term can be estimated, using $(H2)$, we obtain
$$\zeta(t)\int_{0}^{1}(g\circ\theta_{x}){\rm d}x\leq-\int_{0}^{1}(g'\circ\theta_{x}){\rm d}x\leq-\frac{2}{\beta}E'(t).$$
Thus, \eqref{4.3} becomes, for some positive constant $C_{4},$
\bes
\zeta(t)\mathcal{L'}(t)\leq-C\zeta(t)E(t)-C_{4}E'(t).
\lbl{4.4}
\ees
It is clear that
$$L(t)=\zeta(t)\mathcal{L}(t)+C_{4}E(t)\sim E(t).$$
Therefore, using \eqref{4.4} and the fact that $\zeta'(t)\leq0, $ we arrive at
\bes
\L'(t)=\zeta'(t)\mathcal{L}(t)+\zeta(t)\mathcal{L'}(t)+C_{4}E'(t)\leq-C\zeta(t)E(t).
\lbl{4.5}
\ees
A simple integration of \eqref{4.5} over $(t_{0},t) $ leads to
\bes
L(t)\leq L(t_{0})e^{-C\int_{t_0}^{t}\zeta(s){\rm d}s}, \ \ t\geq t_{0}.
\lbl{4.6}
\ees
Recalling \eqref{3.8}, estimate \eqref{4.6} yields the desired result \eqref{2.4}.\quad \quad $\Box$


\end{document}